\documentclass[aop,seceqn,citesort,MSNbibl,dvips]{arximspdf}
\usepackage{mathbh}

\doi{10.1214/08-AOP405}
\volume{37}
\issue{1}
\pubyear{2009}
\firstpage{347}
\lastpage{392}

\makeatletter

\newproclaim{remark}{Remark}
\newtheorem{theorem}{Theorem}[section]
\newtheorem{lemma}[theorem]{Lemma}
\newtheorem{prop}[theorem]{Proposition}

\makeatother

\begin{document}
\begin{frontmatter}

\title{A two cities theorem for the parabolic Anderson~model}
\runtitle{A two cities theorem for the parabolic Anderson model}

\begin{aug}
\author[A]{\fnms{Wolfgang} \snm{K\"onig}\ead[label=e1]{koenig@math.uni-leipzig.de}},
\author[B]{\fnms{Hubert} \snm{Lacoin}\ead[label=e2]{lacoin@math.jussieu.fr}},
\author[C]{\fnms{Peter} \snm{M\"orters}\corref{}\thanksref{t1,t2}\ead[label=e3]{maspm@bath.ac.uk}}\\ and
\author[D]{\fnms{Nadia} \snm{Sidorova}\thanksref{t2}\ead[label=e4]{n.sidorova@ucl.ac.uk}}
\runauthor{W. K\"onig, H. Lacoin, P. M\"orters and N. Sidorova}
\affiliation{Universit\"at Leipzig, Universit\'e Denis Diderot (Paris 7),
University of Bath and University College London}
\address[A]{W. K\"onig\\
Universit\"at Leipzig\\
Fakult\"at f\"ur Mathematik und Informatik\\
Mathematisches Institut\\
Postfach 10 09 20\\
04009 Leipzig\\
Germany\\
\printead{e1}} 
\address[B]{H. Lacoin\\
Universit\'e Denis Diderot (P7)\\
Bo\^\i te courrier 7012\\
75251 Paris Cedex 05\\
France\\
\printead{e2}}
\address[C]{P. M\"orters\\
University of Bath\\
Department of Mathematical Sciences\\
Claverton Down\\
Bath BA2 7AY\\
United Kingdom\\
\printead{e3}}
\address[D]{N. Sidorova\\
University College London\\
Department of Mathematics\hspace*{6.5pt}\\
Gower Street \\
London WC1E 6BT\\
United Kingdom\\
\printead{e4}}
\end{aug}
\thankstext{t1}{Supported by an Advanced Research Fellowship of EPSRC.}
\thankstext{t2}{Supported by Grant EP/C500229/1 of the Engineering and
Physical Sciences Research Council (EPSRC).}
\pdfauthor{Wolfgang Konig, Hubert Lacoin, Nadia Sidorova}

\received{\smonth{10} \syear{2007}}
\revised{\smonth{4} \syear{2008}}

%
\begin{abstract}
The parabolic Anderson problem is the Cauchy problem for the heat equation
$\partial_t u(t,z)=\Delta u(t,z)+\xi(z) u(t,z)$ on $(0,\infty)\times
{\mathbb Z}^d$ with random potential
$(\xi(z) \dvtx z\in{\mathbb Z}^d)$. We consider independent and
identically distributed potentials,
such that the distribution function of $\xi(z)$ converges polynomially
at infinity. If $u$ is initially localized in the origin, that is, if
$u(0,{z})={\mathbh1}_0({z})$,
we show that, as time goes to infinity, the solution is completely localized
in \textit{two} points almost surely and in \textit{one} point with
high probability.
We also identify the asymptotic behavior of the concentration sites in
terms of a weak limit theorem.
\end{abstract}

%
\pdfkeywords{Parabolic Anderson problem,
Anderson Hamiltonian,
random potential,
intermittency,
localization,
pinning effect,
heavy tail,
polynomial tail,
Pareto distribution,
Feynman--Kac formula}

\end{frontmatter}

\section{Introduction and main results}

\subsection{The parabolic Anderson model and intermittency}\label
{sec-intro}

We consider the heat equation with random potential on the integer
lattice $\mathbb Z^d$ and study
the Cauchy problem with localized initial datum,
%
\begin{eqnarray}\label{pam}
\partial_t u(t,z) & = & \Delta u(t,z)+\xi
(z)u(t,z),\qquad
(t,z)\in(0,\infty)\times\mathbb Z^d,\nonumber\\[-8pt]\\[-8pt]
u(0,z) & = & {\mathbh1}_{0}(z), \qquad z\in\mathbb Z^d,\nonumber
\end{eqnarray}
where
\begin{eqnarray*}
(\Delta f)(z)=\sum_{y\sim z}[f(y)-f(z)], \qquad z\in\mathbb Z^d,
f\dvtx
\mathbb Z^d\to\mathbb R
\end{eqnarray*}
is the discrete Laplacian, and the potential
$(\xi(z)\dvtx z\in\mathbb Z^d)$ is a collection
of independent identically distributed random variables.

The problem \textup{(\ref{pam})} and its variants are often called the
\textit{parabolic Anderson problem}. {It}
originated in the work of the physicist P. W. Anderson on entrapment of
electrons
in crystals with impurities, see \cite{An58}. The parabolic version of
the problem appears
in the context of chemical kinetics and population dynamics, and also
provides a simplified
qualitative approach to problems in magnetism and turbulence.
The references \cite{GM90,M94,CM94} provide
applications, background and heuristics around the parabolic Anderson model.
Interesting recent mathematical progress can be found, for example,
in \cite{GH06,GK05,BMR07} and \cite{HKM06} is a recent survey article.

One main reason for the great interest in the parabolic Anderson
problem lies in the fact that it
exhibits an \textit{intermittency effect}: It is believed that, at
late times, the overwhelming contribution to
the total mass of the solution $u$ of the problem~\textup{(\ref
{pam})} comes
from a small number of
{spatially} separated regions of small diameter, which are often called
the \textit{relevant islands}.
As the upper tails of the potential distribution get heavier, this
effect is believed to get stronger and
the number of relevant islands and their sizes are believed to become smaller.
Providing rigorous evidence for
intermittency is a major challenge for mathematicians, which
has lead to substantial research efforts in the past 15 years.

An approach, which has been proposed in the physics literature
(see \cite{GK05} or \cite{ZM87})
suggests that we should study the large time asymptotics of the moments
of the total mass
%
\begin{equation}\label{U(t)def}
U(t)=\sum_{z\in\mathbb Z^d}u(t,z) , \qquad t>0 .
\end{equation}
Denoting expectation with respect to $\xi$ by $\langle\cdot
\rangle$,
if all exponential moments $\langle\exp(\lambda\xi(z)) \rangle$
for $\lambda>0$ exist, then so do
all moments $\langle U(t)^p\rangle$ for $t>0$, $p>0$. Intermittency
becomes manifest in a faster
growth rate of higher moments. More precisely, the model is called
intermittent if
%
\begin{equation}\label{Intermitt}
\limsup_{t\to\infty}\frac{\langle U(t)^p\rangle^{1/p}}{\langle
U(t)^q\rangle^{1/q}}=0\qquad
\mbox{for $0<p<q$.}
\end{equation}
Whenever $\xi$ is nondegenerate random, the parabolic Anderson model
is intermittent in this sense, see \cite{GM90}, Theorem 3.2.
Further properties of the relevant islands, like their asymptotic size and
the shape of potential and solution, are reflected (on a heuristical level)
in the asymptotic expansion of $\log\langle U(t)^p\rangle$
for large $t$. Recently, in \cite{HKM06}, it was argued that the
distributions with finite exponential
moments can be divided into exactly four different universality
classes, with each class having a qualitatively
different long-time behavior of the solution.

It is, however, a much harder mathematical challenge to prove
intermittency in the original
geometric sense, and to identify asymptotically the number, size and
location of the relevant
islands. This program was initiated by Sznitman for the closely related
continuous
model of a Brownian motion with Poissonian obstacles, and the very
substantial body of research he
and his collaborators created is surveyed in his monograph \cite{S98}.
For the problem \textup{(\ref{pam})} and
two universality classes of potentials, the double-exponential
distribution and distributions
with tails heavier than double-exponential (but still with all
exponential moments finite),
the recent paper \cite{GKM07} makes substantial progress toward
completing the geometric picture:
Almost surely, the contribution coming from the complement of a random
number of relevant islands is negligible compared to the mass coming
from these islands, asymptotically
as $t\to\infty$. In the double-exponential case, the radius of the
islands stays bounded; in the heavier
case, the islands are single sites; and in Sznitman's case, the radius
tends to infinity
on the scale $t^{1/(d+2)}$.

Questions about the number of relevant islands remained open in all
these cases, {and constitute
the main concern of the present paper.} In \cite{GKM07,S98} it is
shown that an upper
bound on the number of relevant
islands is $t^{o(1)}$, but this is certainly not always the best
possible bound. In particular, the questions
whether a \textit{bounded number} of islands already carry the bulk of
the mass, or when \textit{just one}
island is sufficient, are unanswered. These questions are difficult,
since there are
many local regions that are good candidates for being a relevant
island, and the known
criterion that identifies relevant islands does not seem to be optimal.

In the present paper, we study the parabolic Anderson model with
potential distributions
that do not have any finite exponential moment.
For such distributions one expects the intermittency effect to be even
more pronounced than in the cases
discussed above, with a very small number of relevant islands, which
are just single sites.
Note that in this case intermittency cannot be studied in terms of the moments
$\langle U(t)^p \rangle$, which are not finite.

The main result of this paper is that, in the case of potentials with
polynomial tails,
\textit{almost surely} at all large times there are \textit{at most
two} relevant islands, each
of which consists of a single site. In other words, the proportion of
the total mass $U(t)$
is asymptotically concentrated in just two time-dependent lattice
points. Note that, by the
intermediate value theorem, the total mass cannot be concentrated in
just one site, if
this site is changing in time on the lattice. Hence this is the
strongest form of localization
that can hold almost surely. However, we also show that, \textit{with
high probability},
the total mass $U(t)$ is concentrated in a \textit{single} lattice point.

The \textit{intuitive picture} is that, at a typical large time, the
mass, which is thought of as a population,
inhabits one site, interpreted as a city. At some rare times, however,
word spreads that a better site has been found,
and the entire population moves to the new site, so that at the
transition times part of the population
still lives in the old city, while another part has already moved to
the new one. This picture inspired the term
\textit{``two cities theorem''} for our main result, which was
suggested to us by S. A. Molchanov. The present
paper is, to the best of our knowledge, the first where such a behavior
is found in a model of
mathematical physics.

Concentration of the mass in a single site with high probability has
been observed so far only for
quite simple mean field models; see \cite{FM,FG}. The present paper
is the first instance where
it has been found in the parabolic Anderson model or, indeed, any
comparable lattice-based model.
We also study the asymptotic locations of the points where the mass
concentrates in terms of
a weak limit theorem with an explicit limiting density. Precise
statements are formulated in the next section.

\subsection{The parabolic Anderson model with Pareto-distributed
potential}\label{sec-result}

We assume that the potentials $\xi(z)$ at all sites $z$ are
independent and \textit{Pareto-distributed}
with parameter $\alpha>d$, that is, the distribution function is
%
\begin{equation}
F(x)=\operatorname{Prob}\bigl(\xi(z)<x\bigr)=1-x^{-\alpha},\qquad x\ge1.
\end{equation}
In particular, we have $\xi(z)\geq1$ for all $z\in\mathbb Z^d$, almost
surely. Note from \cite{GM90}, Theorem 2.1, that the
restriction to parameters $\alpha>d$ is necessary and sufficient
for~\textup{(\ref{pam})} to possess a
unique nonnegative solution $u\dvtx(0,\infty)\times\mathbb Z^d\to
[0,\infty)$. Recall that
\[
U(t)=\sum_{z\in\mathbb Z^d}u(t,z)
\]
is the total mass of the solution at time $t>0$.
Our main result shows the almost sure localization of the solution
$u(t, \cdot)$ in two
lattice points $Z_t^{({1})}$ and $Z_t^{({2})}$, as $t\to\infty$.

\begin{theorem}[(Two cities theorem)]
\label{main_as}
Suppose $u\dvtx(0,\infty)\times\mathbb Z^d\to[0,\infty)$ is the
solution to the parabolic Anderson
problem \textup{(\ref{pam})} with i.i.d. Pareto-distributed
potential with
parameter $\alpha>d$.
Then there exist processes $(Z_t^{({1})}\dvtx t>0)$ and
$(Z_t^{({2})}\dvtx t>0)$
with values in $\mathbb Z^d$, such that
{$Z_t^{({1})}\not= Z_t^{({2})}$ for all $t>0$, and}
\[
\lim_{t\to\infty}\frac{u(t,Z_t^{({1})})+u(t,Z_t^{({2})})}{U(t)}=
1 \qquad\mbox
{almost surely.}
\]
\end{theorem}

\begin{remark}\label{remark1}
At least two sites are
needed to carry the total mass
in an \textit{almost sure} limit theorem. Indeed, assume that
there is a single process $(Z_t \dvtx t>0)$ such that
$u(t,Z_t)>2U(t)/3$ for all large $t$.
As $u( \cdot,z)$ is continuous for any $z\in\mathbb Z^d$, this
leads to a
contradiction at jump times
of the process $(Z_t\dvtx t>0)$. From the growth of $U(t)$ one can see
that this process is not
eventually constant, and thus has jumps at arbitrarily large times.
\end{remark}

Our second result {concerns convergence in probability}. We show that
the solution
$u(t, \cdot)$ is localized in \textit{just one} lattice point
with high probability.

\begin{theorem}[(One point localization in probability)]
\label{main_w}
The process $(Z_t^{({1})}\dvtx t>0)$ in Theorem \ref{main_as} can be chosen
such that
\[
\lim_{t\to\infty}\frac{u(t,Z_t^{({1})})}{U(t)}= 1 \qquad\mbox{in
probability.}
\]
\end{theorem}

\begin{remark}\label{remark2}
The proof of this result given in this paper {uses} strong results
provided for the proof of
Theorem \ref{main_as}. However, it can be proved with less
sophisticated tools, and a
self-contained proof can be found in our unpublished preprint~\cite{KMS06}.
\end{remark}

\begin{remark}\label{remark3}
We conjecture that the one point localization phenomenon holds for a
wider class of
heavy-tailed potentials, including the stretched exponential case. We
also believe that it
does \textit{not} hold for \textit{all} potentials in the
``single-peak'' class of \cite{HKM06}.
\end{remark}

\begin{remark}\label{remark4}
The asymptotic behavior of $\log U(t)$ for the Anderson model with
heavy-tailed potential
is analyzed in detail in \cite{HMS08}. In the case of a
Pareto-distributed potential
it turns out that already the leading term in the asymptotic expansion
of $\log U(t)$ is random.
This is in sharp contrast to potentials with exponential moments, where
the leading two
terms in the expansion are always deterministic. More precisely, introducing
%
\begin{equation}\label{mudef}
q=\frac{d}{\alpha-d} \quad\mbox{and} \quad\theta=\frac{2^d
B(\alpha-d,d)}{{q}^d(d-1)!},
\end{equation}
where $B(\cdot,\cdot)$ denotes the Beta function,
in \cite{HMS08}, Theorem 1.2, it is shown that
%
\begin{equation}\label{U(t)asy}\quad
{\frac{(\log t)^q}{t^{q+1}}} \log U(t)\quad\Longrightarrow\quad
Y\qquad
\mbox{where }
\mathbb P(Y\leq y)=\exp\{-\theta y^{d-\alpha}\}
\end{equation}
%
{and} $\Rightarrow$ denotes weak convergence.
Note that the upper tails of $Y$ have the same asymptotic order as the
Pareto distribution
with parameter $\alpha-d$, that is, $\mathbb P(Y>y)\asymp y^{d-\alpha}$ as
$y\to\infty$.
{The proof} of \cite{HMS08}, Theorem 1.2, also shows that
there is a~process $(Z_t\dvtx t>0)$ such that
%
\begin{equation}\label{u(t,Zt)law}
{\frac{(\log t)^q}{t^{q+1}}} \log u(t,Z_t)\quad\Longrightarrow\quad Y
\qquad\mbox{ where }
\mathbb P(Y\leq y)=\exp\{-\theta y^{d-\alpha}\}.\hspace*{-38pt}
\end{equation}
Note, however, that a combination of \textup{(\ref{U(t)asy})} with
\textup{(\ref{u(t,Zt)law})} does not yield the
concentration property in Theorem \ref{main_w} since the asymptotics
are only
logarithmic. Much more precise techniques are necessary for this purpose.
\end{remark}

In Section \ref{overview} we see how the process $(Z_t^{({1})}\dvtx
t>0)$ in
Theorem \ref{main_w} can be
defined as the maximizer in a random variational problem {associated
with the parabolic Anderson
problem.} Our third result is a limit theorem for this process. Recall
the definition of $q$ and $\theta$
from \textup{(\ref{mudef})}, and denote by $|\cdot|$ the $\ell
^1$-norm on $\mathbb R^d$.

\begin{theorem}[(Limit theorem for {the} concentration site)]
\label{main_z}
The process $(Z_t^{({1})}\dvtx t>0)$ in Theorem \ref{main_w} can be chosen
such that, as $t\to\infty$,
\begin{eqnarray*}
Z_t^{({1})}\biggl(\frac{\log t}{t}\biggr)^{q+1}\quad\Longrightarrow\quad X^{({1})},
\end{eqnarray*}
where $X^{({1})}$ is an $\mathbb R^d$-valued random variable with density
\begin{eqnarray*}
p^{({1})}(x_1)  =\alpha\int_0^{\infty}\frac{\exp\{-\theta
y^{d-\alpha}\}\,dy}
{(y+q|x_1|)^{\alpha+1}}.
\end{eqnarray*}
\end{theorem}

\begin{remark}\label{remark5}
The proof of this result
uses the point process technique
developed in \cite{HMS08}. A more elementary proof can be found in our
unpublished preprint~\cite{KMS06}.
\end{remark}

\begin{remark}\label{remark6}
If we choose the processes $(Z_t^{({1})}\dvtx t>0)$ and
$(Z_t^{({2})}\dvtx t>0)$
such that, with probability tending to one,
$u(t,Z_t^{({1})})$ and $u(t,Z_t^{({2})})$ are the largest and second
largest value of
$u(t,z)$,
we show that, as $t\to\infty$,
\begin{eqnarray*}
\bigl(Z_t^{({1})}, Z_t^{({2})}\bigr) \biggl(\frac{\log t}{t}
\biggr)^{q+1}\quad\Longrightarrow\quad
\bigl(X^{({1})},X^{({2})}\bigr),
\end{eqnarray*}
where $(X^{({1})},X^{({2})})$ is a pair of $\mathbb R^d$-valued random variables
with joint density
\begin{eqnarray*}
p(x_1,x_2)=\int_0^{\infty} \frac{\alpha\exp\{-\theta
y^{d-\alpha}\}\,dy}
{(y+q|x_1|)^{\alpha}(y+q|x_2|)^{\alpha+1}}.
\end{eqnarray*}
By projecting this result on the first component we obtain the
{convergence in distribution statement} 
of Theorem \ref{main_z}, where the density of $X^{({1})}$ is given by
\begin{eqnarray*}
p^{({1})}(x_1)
=\int_0^{\infty}\biggl(\int_{\mathbb R^d} \frac{d
x_2}{(y+q|x_2|)^{\alpha+1}} \biggr)
\frac{\alpha\exp\{-\theta y^{d-\alpha}\}}{(y+q|x_1|)^{\alpha}}\, d
y .
\end{eqnarray*}
The inner integral equals $y^{d-\alpha-1} 2^dq^{-d} B(\alpha+1-d,d)/(d-1)!$.
Recalling \textup{(\ref{mudef})} and using the functional
equation $B(x+1,y) (x+y)=B(x,y) x$ for $x,y>0$, yields
\begin{eqnarray*}
p^{({1})}(x_1) =(\alpha-d)\theta \int_0^{\infty} y^{d-\alpha-1}
\frac{\exp\{-\theta y^{d-\alpha}\}}{(y+q|x_1|)^{\alpha}} \,dy
=\alpha\int_0^{\infty}\frac{\exp\{-\theta y^{d-\alpha}\}\,dy}
{(y+q|x_1|)^{\alpha+1}},
\end{eqnarray*}
using integration by parts in the last step.
Moreover, from the proof of Theorem~\ref{main_z} one can easily infer
the joint
convergence
%
\[
\biggl(\frac{\log t}{t}\biggr)^{q+1} \biggl( Z_t^{({1})}, \frac{\log
u(t,Z_t^{({1})}
)}{\log t} \biggr)
\quad\Longrightarrow\quad(X,Y) ,
\]
where the joint density of $(X,Y)$ is
\[
(x,y)\mapsto\alpha\frac{\exp\{-\theta y^{d-\alpha}\}
}{(y+q|x|)^{\alpha+1}}.
\]
\end{remark}

\subsection{Overview: the strategy behind the proofs}\label
{overview}

{Through{out} the paper we will say that a statement occurs \textit{
eventually for all $t$} when there exists a time $t_0$ such that the
statement is fulfilled for all $t>t_0$. Note that when a statement is
said to {hold} true almost surely eventually for all $t$, {the
corresponding} $t_0$ can be random.}

As shown in \cite{GM90}, Theorem 2.1, under the assumption $\alpha
>d$, the
unique nonnegative solution $u\dvtx(0,\infty)\times\mathbb Z^d\to
[0,\infty)$
of \textup{(\ref{pam})} has a \textit{Feynman--Kac representation}
%
\begin{equation}\label{FKform}
u(t,z)=\mathbb E_0 \biggl[\exp\biggl\{\int_0^t\xi(X_s)\,ds\biggr\}
{\mathbh1}\{
{X_t=z\}}\biggr],\qquad t>0,\ z\in\mathbb Z^d,
\end{equation}
where $(X_s \dvtx s\ge0)$ under $\mathbb P_0$ (with expectation
$\mathbb E_0$) is a
continuous-time simple random walk on the lattice $\mathbb Z^d$ with
generator $\Delta$ starting at the origin.
Hence, the total mass of the solution is given by
\begin{eqnarray*}
U(t)=\mathbb E_0 \biggl[\exp\biggl\{\int_0^t\xi(X_s)\,ds\biggr\}\biggr].
\end{eqnarray*}
Heuristically, for a fixed time $t>0$, the {paths} $(X_s \dvtx0 \le s
\le t)$
that have the greatest impact on the average $U(t)$ {spend most of
their time at a site $z$}, which has a large potential value $\xi(z)$
and can be reached quickly, that is, is sufficiently close to the origin.

For $\rho\in(0,1)$, the strategy $A^{z,\rho}_t$ of wandering to a
site $z$ during the time interval $[0,\rho t)$
and staying at $z$ during the time $[\rho t,t]$ has, for $|z|\gg t$,
approximately the probability
\begin{eqnarray*}
\mathbb P_0(A^{z,\rho}_t)\approx\exp\biggl\{-|z|\log
\frac
{|z|}{e\rho t}+\eta(z)\biggr\},
\end{eqnarray*}
where $\eta(z)=\log N(z)$ and $N(z)$ denotes the number of paths of
length $|z|$ starting
at zero and ending at $z$ {(see Proposition \ref{iub} for details)}.
Then the integral in the Feynman--Kac formula is bounded from below
by $t(1-\rho)\xi(z)$ for the paths of the random walk following the
strategy $A^{z,\rho}_t$. Hence,
we obtain by optimizing over $z$ and $\rho\in(0,1)$,
\begin{eqnarray*}
\frac{1}{t}\log U(t)
&\gtrapprox&\sup_{z\in\mathbb Z^d}
\sup_{\rho\in(0,1)}\biggl[(1-\rho)\xi(z)-\frac{|z|}{t}\log\frac
{|z|}{e\rho t}+\frac{\eta(z)}{t}\biggr]
\\
&=&\max_{z\in\mathbb Z^d}\Phi_t(z),
\end{eqnarray*}
where
%
\begin{eqnarray}
\label{phi}
\Phi_t(z)=\biggl[\xi(z)-\frac{|z|}{t}\log\xi(z)+\frac{\eta
(z)}{t}\biggr]{\mathbh1}\{t\xi(z)\ge|z|\}.
\end{eqnarray}
The restriction $t\xi(z)\ge|z|$ arises {as, otherwise, the globally
optimal value
$\rho=|z|/(t\xi(z))$ would exceed one.}
This bound, {stated as Proposition \ref{iub}}, is a minor improvement
of the {lower} bound obtained in \cite{HMS08}.
{In addition,} we show that $\max\Phi_t$ also gives an asymptotic
upper bound
for $\frac1t \log U(t)$, which is much harder and constitutes a
significant improvement of the bound obtained in \cite{HMS08}; see
Proposition \ref{cor-uppboundU}. Altogether
\begin{eqnarray*}
\frac{1}{t}\log U(t)\approx\max_{z\in\mathbb Z^d}\Phi_t(z)
\end{eqnarray*}
and it is plausible that the optimal sites at time $t$ are the sites
where the two largest
values of the random functional $\Phi_t$ are attained.
This is indeed the definition of the processes $(Z_t^{({1})}\dvtx t\ge
0)$ and
$(Z_t^{({2})}\dvtx t\ge0)$, which underlies our three main theorems.

The remainder of the paper is organized as follows.

In Section \ref{s_pre} we provide several technical results for later use.
In particular, we study the behavior of $\eta(z)$ and of the upper
order statistics
of the potential $\xi$, and we derive spectral estimates similar to
those obtained in \cite{GKM07}.

In Section \ref{s_phi} { we study the asymptotic properties of the
sites $Z_t^{({1})}$, $Z_t^{({2})}$ and $Z_t^{({3})}$, where $\Phi_t$
attains its three
largest values, as well as the properties
of $\Phi_t(Z_t^{({i})})$, for $i=1,2,3$}. Here we prove
Proposition \ref{gap}, which
states that, almost surely, the gap $\Phi_t(Z_t^{({1})})-\Phi
_t(Z_t^{({3})})$ is eventually large
enough. This is the main reason
for $u(t,z)$ being concentrated at just two sites $Z_t^{({1})}$ and
$Z_t^{({2})}$.
Observe that a similar statement
about the gap $\Phi_t(Z_t^{({1})})-\Phi_t(Z_t^{({2})})$ is not true
as, {by continuity}, there are
arbitrarily large times $t$
such that $\Phi_t(Z_t^{({1})})=\Phi_t(Z_t^{({2})})$, which is the
main technical reason for the
absence of one point almost sure localization.

In Section \ref{s_mass} we study the total mass of the solution and
its relation to $\Phi_t$.
We split $U(t)$ into five parts according to five groups of paths, and
show that only one of
them {makes an essential contribution, namely the one corresponding} to
paths which
visit either $Z_t^{({1})}$ or $Z_t^{({2})}$ and whose length is not too
large. Then
we prove Propositions \ref{iub}
and \ref{cor-uppboundU}, which are the very precise upper and lower
approximations
of $\frac{1}{t}\log U(t)$ by $\Phi_t(Z_t^{({1})})$ needed for
Theorem \ref{main_as}.

In Section \ref{s_as} we prove Theorem \ref{main_as}. We split the
probability space into three
disjoint events:
\begin{itemize}
\item The gap $\Phi_t(Z_t^{({1})})-\Phi_t(Z_t^{({2})})$ is small and
the sites $Z_t^{({1})}$ and $Z_t^{({2})}$ are close.
\item The gap $\Phi_t(Z_t^{({1})})-\Phi_t(Z_t^{({2})})$ is small but
the sites $Z_t^{({1})}$ and $Z_t^{({2})}$ are
far away.
\item The gap $\Phi_t(Z_t^{({1})})-\Phi_t(Z_t^{({2})})$ is large.
\end{itemize}
Correspondingly, we prove Propositions \ref{ttwo}, \ref{two} and \ref
{three}, which justify
Theorem~\ref{main_as} for each event. In each case, we decompose
$u(t,z)$ in two components
(differently for different events) and show that one of them localizes
around $Z_t^{({1})}$
and $Z_t^{({2})}$, and the other one is negligible.

Finally, in Section \ref{s_weak} we prove Theorems \ref{main_w}
and \ref{main_z}.
We use the point processes technique developed in \cite{HMS08}, which
readily gives
Theorem \ref{main_z}. Theorem \ref{main_w} is obtained using a
combination of the point processes
approach and Theorem \ref{main_as}.

\section{Notation and preliminary results}\label{s_pre}

For $z\in\mathbb Z^d$, we define $N(n,z)$ as the number of paths of length
$n$ in $\mathbb Z^d$ starting at the
origin and passing through $z$. Recall that $N(z)=N(|z|,z)$, {where
here and throughout
the paper $|\cdot|$ denotes the $\ell^1$-norm.}
For $n\ge|z|$, we define
\begin{eqnarray*}
\eta(n,z)=\log N(n,z) \quad\mbox{and}\quad \eta(z)=\log N(z).
\end{eqnarray*}
It is easy to see that $0\le\eta(z)\le|z|\log d$. We define two
important scaling functions
%
\begin{equation}\label{rtatdef}
r_t=\biggl(\frac{t}{\log t}\biggr)^{q+1}
\quad\mbox{and}\quad
{a}_t=\biggl(\frac{t}{\log t}\biggr)^{q},
\end{equation}
where $r_t$ will turn out to be the appropriate scaling for $Z_t^{({i})}$
and $a_t$ for $\Phi_t(Z_t^{({i})})$, $i=1,2,3$.

For each $r>0$, denote $\xi_r^{({1})}=\max_{|z|\le r}\xi
(z)$ and
\begin{eqnarray*}
\xi^{({i})}_r=\max\bigl\{\xi(z)\dvtx|z|\le r, \xi(z)\neq\xi
_r^{({j})}\ \forall
{j< i} \bigr\}
\end{eqnarray*}
for $2\le i\le\ell_r$, where $\ell_r$ is the number of points {in
the ball $\{|z|\le r\}$.}
Hence,
\[
\xi_r^{({1})}>\xi_r^{({2})}>\xi_r^{({3})}>\cdots>\xi_r^{({\ell_r})}
\]
are precisely the potential values in {this} ball.

Fix $0<\rho<\sigma<\frac{1}{2}$ so that $\sigma<1-\frac{\rho
}{d}$, and
$\nu>0$. We define four auxiliary scaling {functions}
\begin{eqnarray*}
f_t&=&(\log t)^{-{1}/{d}-\nu},\qquad g_t=(\log t)^{{1}/({\alpha
-d})+\nu},\\
k_t&=&\lfloor{\lfloor} r_tg_t{\rfloor}^{\rho} \rfloor,\qquad\hspace*{13.75pt} m_t=
\lfloor{\lfloor} r_tg_t{\rfloor}^{\sigma} \rfloor
\end{eqnarray*}
and two sets
\begin{eqnarray*}
F_t&=&\bigl\{z\in\mathbb Z^d\dvtx|z|\le r_tg_t, \exists i<k_t \mbox{ such
that }\xi(z)=\xi_{r_tg_t}^{({i})}\bigr\},\\
G_t&=&\bigl\{z\in\mathbb Z^d\dvtx|z|\le r_tg_t, \exists i<m_t \mbox{ such
that }\xi(z)=\xi_{r_tg_t}^{({i})}\bigr\},
\end{eqnarray*}
which will be used throughout this paper. In other words, $F_t$,
respectively, $G_t$, is the set of those sites in
{the ball $\{|z|\le r_tg_t\}$} in which the $k_t-1$, respectively,
$m_t-1$, largest potential sites
are attained. Hence $F_t\subset G_t$ and $F_t$, respectively, $G_t$,
have precisely $k_t-1$, respectively, $m_t-1$,
elements.

\subsection{Two technical lemmas}

We start by proving an estimate on $\eta(n,z)$, which
we will use later in order to prove that if $z$ is a point where the
potential is high,
then a path passing through $z$ only contributes to the Feynman--Kac
formula if its length is close to $|z|$.

\begin{lemma}
\label{nn}
There is a constant $K$ such that for all $n\ge|z|$,
\begin{eqnarray*}
\eta(n,z)-\eta(z)\le(n-|z|)\log\frac{2den}{n-|z|}+K.
\end{eqnarray*}
\end{lemma}

\begin{pf} We fix $z=(z_1,\ldots,z_d)\in\mathbb Z^d$ and without loss
of generality
assume that $z_i\ge0$. Denote by $\mathcal{P}_{n,z}$ the set of paths
of length $n$
starting at the origin and passing through $z$. Each $y\in\mathcal{P}_{n,z}$
can be described by the vector $(y_1,\ldots,y_n)$ of its increments,
where $|y_i|=1$ for all $i$.
Since the path $y$ passes through $z$, there is a subsequence
$(y_{i_1},\ldots,y_{i_{|z|}})$
corresponding to a path from $\mathcal{P}_{|z|,z}$. Thus, every path
from $\mathcal{P}_{n,z}$
can be obtained from a path in $\mathcal{P}_{|z|,z}$ by adding $n-|z|$
elements to its
coding sequence. As there are only $2d$ possible elements and
${n\choose n-|z|}$
possibilities where the elements can be added, we obtain an upper bound
\begin{eqnarray*}\hspace*{3pt}
N(n,z)\le N(z)(2d)^{n-|z|}\pmatrix{n\cr n-|z|}\le N(z)\frac
{(2dn)^{n-|z|}}{(n-|z|)!}
\le N(z)e^K\biggl(\frac{2den}{n-|z|}\biggr)^{n-|z|}
\end{eqnarray*}
with $K$ such that $m!\ge e^{-K}(m/e)^m$ for all $m$.
Taking the logarithm completes the proof.
\end{pf}

In the next\vspace*{1pt} lemma we derive some properties of the upper order statistics
of the potential $\xi$, which will be used later to prove that $\Phi
_t(Z_t^{({1})})$
is an approximate
upper bound for
$\frac{1}{t}\log U(t)$.

\begin{lemma}
\label{asymp}
There exists $c>0$ such that, with probability one, eventually for all
$t$:
\begin{longlist}
\item  $t^c<\xi_{r_tg_t}^{({k_t})}<t^{q-c}$
and $\xi_{r_tg_t}^{({m_t})}/\xi_{r_tg_t}^{({k_t})}<t^{-c}$;
\item $F_t\cap\{|z|\le t^{(q+1)\sigma+c}\}=\varnothing$;
\item  $G_t$ is totally disconnected, that is, if $x,y\in
G_t$, then $|x-y|\neq1$.
\end{longlist}
\end{lemma}

\begin{pf} (i) Note that $\widehat\xi(z)=\alpha\log\xi(z)$
defines a field of independent exponentially
distributed random variables. It has been proved in \cite
{HMS08}, (4.7), that, for each $\kappa\in(0,1)$,
\begin{eqnarray*}
\lim_{n\to\infty}\frac{\log\xi_n^{({\lfloor n^\kappa\rfloor
})}}{\log n}=\frac{d-\kappa}{\alpha}\qquad \mbox{almost surely.}
\end{eqnarray*}
Substituting $n=r_tg_t$ and $\kappa=\rho$, respectively, $\kappa
=\sigma$, we obtain
%
\begin{eqnarray}\label{km}
\lim_{t\to\infty}\frac{\log\xi_{r_tg_t}^{({k_t})}}{\log t}
&=&\frac{(d-\rho)(q+1)}{\alpha}
\quad\mbox{and}\nonumber\\[-8pt]\\[-8pt]
\lim_{t\to\infty}\frac{\log\xi
_{r_tg_t}^{({m_t})}}{\log t}&=&\frac{(d-\sigma)(q+1)}{\alpha}.\nonumber
\end{eqnarray}
{The result follows, since $\frac{(d-\rho)(q+1)}{\alpha}\in(0,q)$
for $\rho\in(0,1)$ and
$\frac{(d-\rho)(q+1)}{\alpha}>\frac{(d-\sigma)(q+1)}{\alpha}$
for $\rho<\sigma$.}%

\hspace*{3pt}(ii) Because $\sigma<1-\frac{\rho}{d}$, we can pick $c$ and
$\varepsilon$ small enough such that
$\sigma+\frac{cd}{q\alpha}+\frac{2\varepsilon}{q}<1-\frac{\rho
}{d}$. Then by \cite{HMS08}, Lemma 3.5, we obtain
\begin{eqnarray*}
\max_{|z|\le t^{(q+1)\sigma+c}}\xi(z)\le t^{{d}/{\alpha
}[(q+1)\sigma+c]+\varepsilon}
=t^{q\sigma+{cd}/{\alpha}+\varepsilon}
<t^{{(d-\rho)q}/{d}-\varepsilon}
\end{eqnarray*}
eventually, which, together with {the first part of \textup{(\ref{km})}}
implies the statement.

(iii) For each $n\in\mathbb N$, denote $h_n=\lfloor n^{\sigma
}\rfloor$ and
\begin{eqnarray*}
\widehat G_n=\bigl\{z\in\mathbb Z^d\dvtx|z|\le n, \exists i<h_n
\mbox{ such that }\xi(z)=\xi_n^{({i})}\bigr\}.
\end{eqnarray*}
Since $G_t=\widehat G_{\lfloor r_tg_t\rfloor}$, it suffices to show
that $\widehat G_n$
is totally disconnected eventually.

First, consider {the case} $d\ge2$. The set $\widehat G_n$ consists of $h_n$
different points belonging to the ball {$B_n=\{|z|\le n\}$}. Denote
them by
$a_0,\ldots,a_{h_n-1}$, where $a_i$ is such that $\xi(a_i)=\xi
_n^{({i})}$.
For $i\neq j$, the pair $(a_i,a_j)$ is uniformly distributed over all
pairs of
distinct points in $B_n$. Hence the probability of $a_i$ and $a_j$
being neighbors, written $a_i\sim a_j$, can be estimated by
\begin{eqnarray*}
\operatorname{Prob}(a_i\sim a_j)\le\max_{|z|\le
n}\operatorname{Prob}(a_i\sim z
| a_j=z)
\le\frac{2d}{\ell_n-1},
\end{eqnarray*}
where $\ell_n$ is the number of points in $B_n$.
Summing over all pairs, we get
%
\begin{eqnarray}
\label{sim}
\operatorname{Prob}(\widehat G_n\mbox{ is not totally
disconnected})
&\le& \sum_{0\le i<j<h_n} \operatorname{Prob}
(a_i\sim a_j)\nonumber\\[-8pt]\\[-8pt]
&\le&\frac{2dh_n^2}{\ell_n-1}\le Cn^{2\sigma-d}\nonumber
\end{eqnarray}
for some $C>0$. As $\sigma<1/2$ and $d\ge2$, this sequence is summable.
By the Borel--Cantelli lemma $\widehat G_n$ is eventually totally disconnected.

The situation is more delicate if $d=1$. Pick $\sigma'\in(\sigma,1/2)$
and denote $h'_n=\lfloor n^{\sigma'}\rfloor$ and
\begin{eqnarray*}
\widehat G'_n=\bigl\{z\in\mathbb Z^d\dvtx|z|\le n, \exists i<h'_n
\mbox{ such that }\xi(z)=\xi_n^{({i})}\bigr\}.
\end{eqnarray*}
Further, let $p_n=2^{\lfloor\log_2 n\rfloor}$ 
such that $p_n\le n<2p_n$.

It is easy to see
that $\widehat G'_{p_n}$ is totally disconnected eventually.
Indeed, \textup{(\ref{sim})}
remains true with ${\hat G_n}$ and $h_n$ replaced by
${\hat G'_n}$ and $h'_n$,
respectively, and one just needs to notice that $\sum_{n=1}^{\infty
}2^{n(2\sigma'-d)}<\infty$
for $d=1$ and $\sigma'<1/2$.

The final step is to prove that $\widehat G_n\subset\widehat
G'_{2p_n}$. Let
{$\varkappa_n$ be the cardinality of $\widehat G'_{2p_n}\cap B_n$ and
observe that,
for this purpose, it suffices to show that $\varkappa_n\ge h_n$.
Indeed, on this set the $\varkappa_n$ largest values
of $\xi$ over $B_n$ are achieved. We actually prove} a stronger
statement, showing that there are at least $h_{2p_n}$ points
from $\widehat G'_{2p_n}$ in the ball $B_{p_n}$. {From now on} we drop
the subscript $n$.
We write
\[
\widehat G'_{2p}=\{a'_0,\ldots,a'_{h'_{2p}-1}\},
\]
where {$a'_i$} is such that
$\xi(a'_i)=\xi^{({i})}_{2p}$. Let $X=(X_i\dvtx0\le i<h'_{2p})$
with $X_i={\mathbh1}\{|a'_i|\le p\}$ and
\begin{eqnarray*}
|X|=\sum_{i=0}^{h'_{2p}-1}X_i.
\end{eqnarray*}
Since $h'_{2p}=o(p)$ and $|B_p|=2p+1$, $|B_{2p}|=4p+1$, we obtain,
using that the points in $\widehat G'_{2p}$ are uniformly distributed
over $B_{2p}$ without repetitions, that
for large~$p$
\begin{eqnarray*}
\operatorname{Prob}(X_j=1\mid X_i=x_i\ \forall i<j)<3/4
\end{eqnarray*}
and
\begin{eqnarray*}
\operatorname{Prob}(X_j=0\mid X_i=x_i\ \forall i<j)<3/4
\end{eqnarray*}
for all $j<h'_{2p}$ and all $(x_0,\ldots,x_{j-1})\in\{0,1\}^{j}$.
Hence, for all $x\in\{0,1\}^{h'_{2p}}$,
\begin{eqnarray*}
\operatorname{Prob}(X=x)\le(3/4)^{h'_{2p}}.
\end{eqnarray*}
This yields
\begin{eqnarray*}
\operatorname{Prob}(|X|<h_{2p})
&=&\sum_{i=0}^{h_{2p}-1}\sum_{|x|=i}\operatorname{Prob}
(X=x)\\
&\le&\sum_{i=0}^{h_{2p}-1}\pmatrix{h'_{2p}\cr i}(3/4)^{h'_{2p}}
\le h_{2p}(h'_{2p})^{h_{2p}-1}(3/4)^{h'_{2p}}\\
&\le&\exp\{-h'_{2p}\log(4/3)+h_{2p}\log h'_{2p}+\log
h_{2p}\}
=e^{-c(2p)^{\sigma'}}
\end{eqnarray*}
for some $c>0$. Since this sequence is summable, 
we have $|X|\ge h_{2p}$ even-\break tually.
\end{pf}

\subsection{Spectral estimates}\label{sec-spectral}

In this section we exploit ideas developed in \cite{GKM07}.
Let $A\subset\mathbb Z^d$ be a bounded set and denote by $Z_A\in A$
the point,
where the potential $\xi$ takes its maximal value over $A$. Denote by
\begin{eqnarray*}
\mathfrak{g}_A=\xi(Z_A)-\max_{z\in A\setminus\{Z_A\}}\xi(z)
\end{eqnarray*}
the gap between the largest value and the rest of the potential on $A$.
Denote by~$A^*$ the connected component of $A$ containing $Z_A$.
Let $\gamma_A$ and $v_A$ be the principal eigenvalue and eigenfunction
of $\Delta+\xi$ with zero boundary conditions in $A^{*}$ extended by
zero to the whole
set $A$. We assume that $v_A$
is normalized to $v_A(Z_A)=1$. Recall that { under $\mathbb P_z$ and
$\mathbb E_z$
the process
$(X_t\dvtx t\in[0,\infty))$ is a simple random walk with generator
$\Delta$ started from $z\in\mathbb Z^d$.}
The entrance time to a set $A$ is denoted $\tau_A=\inf\{t\geq0\dvtx
X_t\in A\}$, and we write
$\tau_z$ instead of $\tau_{\{z\}}$. Then, {as in \cite{GKM07}}, (4.4),
the eigenfunction $v_A$ admits the probabilistic representation
%
\begin{eqnarray}
\label{ef}\quad
v_A(z)=\mathbb E_z\biggl[\exp\biggl\{\int_0^{\tau_{Z_A}}[\xi
(X_s)-\gamma
_A]\,ds\biggr\}
{\mathbh1}\{\tau_{Z_A}<\tau_{A^{\mathrm{c}}}\}\biggr],\qquad z\in A.
\end{eqnarray}
It turns out that $v_A$ is concentrated around the maximal point $Z_A$
of the potential.

\begin{lemma}
\label{spectral0}
There is a decreasing function $\varphi\dvtx(2d,\infty)\to\mathbb R_{+}$
such that
$\lim_{x\to\infty}\varphi(x)=0$ and,
for any bounded set $A\subset\mathbb Z^d$ satisfying $\mathfrak{g}_A>2d$,
\begin{eqnarray*}
\Vert v_A\Vert_2^2\sum_{z\in A\setminus\{Z_A\}} v_A(z)\le\varphi
(\mathfrak{g}_A).
\end{eqnarray*}
\end{lemma}

\begin{pf} It suffices to consider $z\in A^*$. By the Rayleigh--Ritz
formula we have
\begin{eqnarray*}
\gamma_A&=&\sup\{\langle(\Delta+\xi)f,f\rangle\dvtx f\in
\ell^2(\mathbb Z^d), \operatorname{supp}(f)\subset A^*,
\|f\|_2=1\}\\
&\ge&\sup\{\langle(\Delta+\xi)\delta_z,\delta_z\rangle
\dvtx z\in A^*\}
=\sup\{\xi(z)-2d\dvtx z\in A^*\}\\
&=&\xi(Z_A)-2d.
\end{eqnarray*}
Since the paths of the random walk $(X_s)$ in \textup{(\ref{ef})} do
not leave
$A$ and avoid the point~$Z_A$ where the maximum of the potential is achieved, we can estimate
the integrand
using the gap $\mathfrak{g}_A$. Hence, we obtain
\begin{eqnarray*}
v_A(z)\le\mathbb E_z\bigl[\exp\bigl\{\tau_{Z_A}\bigl(\xi(Z_A)-\mathfrak
{g}_A-\gamma_A\bigr)\bigr\}\bigr]
\le\mathbb E_z[\exp\{-\tau_{Z_A}(\mathfrak{g}_A-2d)\}
].
\end{eqnarray*}
Under $\mathbb P_z$ the random variable $\tau_{Z_A}$ is stochastically
bounded from below by a sum of
$|z-Z_A|$ independent exponentially distributed random variables with
parameter~$2d$. If $\tau$
denotes such a random time, we therefore have
\begin{eqnarray*}
v_A(z)\le\bigl(\mathbb E\bigl[e^{-\tau(\mathfrak
{g}_A-2d)}\bigr]
\bigr)^{|z-Z_A|}=\biggl(\frac{2d}{\mathfrak{g}_A}\biggr)^{|z-Z_A|}.
\end{eqnarray*}
The statement of the lemma follows easily with
\begin{eqnarray*}
\varphi(x)=\Biggl(\sum_{z\in\mathbb Z^d}(2d/x)^{2|z|}\Biggr)\Biggl(\sum
_{z\in
\mathbb Z^d\setminus\{0\}}(2d/x)^{|z|}\Biggr),
\end{eqnarray*}
which obviously satisfies the required conditions.
\end{pf}

Let now $B\subset\mathbb Z^d$ be a bounded set containing the origin and
$\Omega\subset B$. Denote
\begin{eqnarray*}
\mathfrak{g}_{\Omega,B}=\min_{z\in\Omega}\xi(z)-\max_{z\in
B\setminus\Omega}\xi(z)
\end{eqnarray*}
and denote, {for any $(t,z)\in(0,\infty)\times\mathbb Z^d$,}
\begin{eqnarray*}
u_{\Omega,B}(t,z)=\mathbb E_0\biggl[\exp\biggl\{\int_0^t\xi
(X_s)\,ds\biggr\}
{\mathbh1}\{X_t=z\} {\mathbh1}\{\tau_{\Omega}\le t,\tau
_{B^{\mathrm{c}}}>t\}
\biggr].
\end{eqnarray*}

\begin{lemma}
\label{spectral}
Assume that $\mathfrak{g}_{\Omega,B}>2d$. Then, for all $z\in\mathbb Z^d$
and $t>0$:
\begin{itemize}
\item[(a)]
$u_{\Omega,B}(t,z)\le\sum_{y\in\Omega}u_{\Omega,B}(t,y)
\Vert v_{(B\setminus\Omega)\cup\{y\}}\Vert_2^2
v_{(B\setminus\Omega)\cup\{y\}}(z),$
\item[{(b)}]
$\frac{\sum_{z\in B\setminus\Omega}u_{\Omega,B}(t,z)}{\sum_{z\in
B}u_{\Omega,B}(t,z)}\le\varphi(\mathfrak{g}_{\Omega,B}).$
\end{itemize}
\end{lemma}

\begin{pf} (a) This is a slight generalization of \cite{GKM07}, Theorem
4.1, with a ball replaced by an arbitrary bounded set $B$; we repeat the
proof here for the sake of completeness.
For each $y\in\Omega$, by time reversal and using the Markov property
at time $s$, we obtain a lower bound for $u_{\Omega,B}(t,y)$
by requiring that the random walk (now started at $y$) is at $y$ at
time $u$ and has not entered
$\Omega\setminus\{y\}$ before. We have
%
\begin{eqnarray}\label{sp2}\qquad
u_{\Omega,B}(t,y)
&=&\mathbb E_y\biggl[\exp\biggl\{\int_0^t\xi(X_s)\,ds\biggr\}
{\mathbh1}\{X_t=0\}{\mathbh1}\bigl\{\tau_{\Omega\setminus\{y\}}\le
t,\tau
_{B^{\mathrm{c}}}>t\bigr\}\biggr]\nonumber\\
&\ge&\mathbb E_y\biggl[\exp\biggl\{\int_0^u\xi(X_s)\,ds\biggr\}
{\mathbh1}\{X_u=y\}{\mathbh1}\bigl\{\tau_{\Omega\setminus\{y\}}>
u,\tau_{B^{\mathrm
{c}}}>u\bigr\}\biggr]\\
&&{} \times\mathbb E_y\biggl[\exp\biggl\{\int_0^{t-u} \xi
(X_s)\,ds\biggr\}
{\mathbh1}\{X_{t-u}=0\}{\mathbh1}\{\tau_{B^{\mathrm{c}}}>t-u\}\biggr].\nonumber
\end{eqnarray}
Using an eigenvalue expansion for the parabolic problem in
$(B\setminus\Omega)\cup\{y\}$
represented by the first factor on the right-hand side of the formula
above, we obtain the bound
\begin{eqnarray*}
&&\mathbb E_y\biggl[\exp\biggl\{\int_0^u\xi(X_s)\,ds\biggr\}
{\mathbh1}\{X_u=y\}{\mathbh1}\bigl\{\tau_{\Omega\setminus\{y\}}>
u,\tau_{B^{\mathrm
{c}}}>u\bigr\}\biggr]\\
&&\qquad\ge e^{u\gamma_{(B\setminus\Omega)\cup\{y\}}}
\frac{v_{(B\setminus\Omega)\cup\{y\}}(y)^2}{\Vert v_{(B\setminus
\Omega)\cup\{y\}}\Vert_2^2},
\end{eqnarray*}
where we have used that $Z_{(B\setminus\Omega)\cup\{y\}}=y$ since
$\mathfrak{g}_{\Omega,B}>0$.
Substituting the above estimate into \textup{(\ref{sp2})} and taking into
account that
$v_{(B\setminus\Omega)\cup\{y\}}(y)=1$, we obtain
\begin{eqnarray*}
&&\mathbb E_y\biggl[\exp\biggl\{\int_0^{t-u} \xi(X_s)\,ds\biggr\}
{\mathbh1}\{X_{t-u}=0\}{\mathbh1}\{\tau_{B^{\mathrm{c}}}>t-u\}
\biggr]\\
&&\qquad\le
e^{-u\gamma_{(B\setminus\Omega)\cup\{y\}
}}\bigl\Vert v_{(B\setminus
\Omega)\cup\{y\}}\bigr\Vert_2^2u_{\Omega,B}(t,y).
\end{eqnarray*}
The claimed \mbox{estimate} is obvious for $z\notin B$. For $z\in\Omega$,
it follows from\break
$v_{(B\setminus\Omega)\cup\{z\}}(z)=1$, which is implied by
$\mathfrak{g}_{\Omega,B}>0$ and hence
$Z_{(B\setminus\Omega)\cup\{z\}}=z$. Let us now assume that $z\in
B\setminus\Omega$.
Using time reversal, the strong Markov property at time $\tau_{\Omega}$,
and the previous lower bound with $u=\tau_{y}$ we obtain
\begin{eqnarray*}
&&u_{\Omega,B}(t,z)
\\
&&\qquad=\sum_{y\in\Omega}\mathbb E_z\biggl[\exp\biggl\{\int_0^{\tau_y}\xi
(X_s)\,ds\biggr\}
{\mathbh1}\{\tau_y=\tau_{\Omega}\le t,\tau_{B^{\mathrm{c}}}>\tau
_y\}\\
&&\qquad\quad\hspace*{31.7pt}{}\times \mathbb E_y\biggl[\exp\biggl\{\int_0^{t-u}
\xi
(X_s)\,ds\biggr\}
{\mathbh1}\{X_{t-u}=0\}{\mathbh1}\{\tau_{B^{\mathrm{c}}}>t-u\}
\biggr]_{u=\tau
_y}\biggr]\\
&&\qquad\le\sum_{y\in\Omega}u_{\Omega,B}(t,y)\bigl\Vert v_{(B\setminus\Omega
)\cup\{y\}}\bigr\Vert_2^2\\
&&\qquad\quad\hspace*{14.2pt}{}\times
\mathbb E_z\biggl[\exp\biggl\{\int_0^{\tau_y}\bigl[\xi(X_s)-\gamma
_{(B\setminus\Omega)\cup\{y\}}\bigr]\,ds\biggr\}
{\mathbh1}\{\tau_y<\tau_{B^{\mathrm{c}}}\}\biggr]\\
&&\qquad=\sum_{y\in\Omega}u_{\Omega,B}(t,y)\bigl\Vert v_{(B\setminus\Omega)\cup
\{y\}}\bigr\Vert_2^2
v_{(B\setminus\Omega)\cup\{y\}}(z).
\end{eqnarray*}

(b) It suffices to apply Lemma \ref{spectral0} to $A=(B\setminus
\Omega)\cup\{y\}$,
note that $\mathfrak{g}_A\ge\mathfrak{g}_{\Omega,B}$, and use the
monotonicity of $\varphi$. Using
(a), we obtain
\begin{eqnarray*}
\sum_{z\in B\setminus\Omega}u_{\Omega,B}(t,z)
&\le&\sum_{y\in\Omega}u_{\Omega,B}(t,y)\sum_{z\in B\setminus
\Omega}\bigl\Vert v_{(B\setminus\Omega)\cup\{y\}}\bigr\Vert_2^2
v_{(B\setminus\Omega)\cup\{y\}}(z)\\
&\le&\sum_{y\in\Omega}u_{\Omega,B}(t,y) \varphi\bigl(\mathfrak
{g}_{(B\setminus\Omega)\cup\{y\}}\bigr)
\le\varphi(\mathfrak{g}_{\Omega,B})\sum_{y\in B}u_{\Omega,B}(t,y),
\end{eqnarray*}
which completes the proof.
\end{pf}

\section{Properties of the maximizers $Z_t^{({i})}$ and values $\Phi
_t(Z_t^{({i})})$}\label{s_phi}

In this section we introduce the three maximizers $Z_t^{({1})}$,
$Z_t^{({2})}$ and
$Z_t^{({3})}$ and analyze some
{of their crucial properties.} 
In Section \ref{sec-Z1Z2Z3} we concentrate on the {long-term behavior
of the maximizers themselves}
and in Section \ref{sec-Phi1-Phi3} we prove that the maximal value
$\Phi_t(Z_t^{({1})})$ is well separated
from $\Phi_t(Z_t^{({3})})$.

\subsection{The maximizers $Z_t^{({1})}$, $Z_t^{({2})}$
and $Z_t^{({3})}$}\label{sec-Z1Z2Z3}

Recall that $Z_t^{({1})}$, $Z_t^{({2})}$ and $Z_t^{({3})}$ denote the
first three maximizers of
the random functional $\Phi_t$ defined in \textup{(\ref{phi})}. More
precisely,
we define $Z_t^{({i})}$ to be such that
%
\begin{eqnarray}
\label{max}
\Phi_t\bigl(Z_t^{({1})}\bigr)&=&\max_{z\in\mathbb Z^d}\Phi_t(z),\qquad
\Phi_t\bigl(Z_t^{({2})}\bigr)=\max_{z\in\mathbb Z^d\setminus\{Z_t^{({1})}\}}
\Phi_t(z),\nonumber\\[-8pt]\\[-8pt]
\Phi_t\bigl(Z_t^{({3})}\bigr)&=& \max_{z\in\mathbb Z^d\setminus\{Z_t^{({1})}
,Z_t^{({2})}\}} \Phi_t(z).\nonumber
\end{eqnarray}
%

\begin{lemma} With probability one, $Z_t^{({1})}$, $Z_t^{({2})}$ and
$Z_t^{({3})}$ are well
defined for any $t>0$.
\end{lemma}

\begin{pf} Fix $t>0$.
Let $\varepsilon\in(0,1-\frac d\alpha)$. By \cite{HMS08}, Lemma 3.5, there
exists a~random radius $\rho(t)>0$ such that, almost surely,
%
\begin{eqnarray}
\label{asas}
\xi(z)\le\xi_{|z|}^{({1})}\le|z|^{{d}/{\alpha}+\varepsilon
}\le\frac{|z|}{t} \qquad\mbox{for all }|z|>\rho(t).
\end{eqnarray}
Consider $|z|>\max\{\rho(t),edt\}$.
If $t\xi(z)<|z|$ then $\Phi_t(z)=0$. Otherwise, using $\eta(z)\le
|z|\log d$
and estimating $\xi(z)$ in two different ways, we obtain
\begin{eqnarray*}
\Phi_t(z)\le\xi_{|z|}^{({1})}-\frac{|z|}{t}\log\frac
{|z|}{dt}\le
\frac{|z|}{t}\biggl[1-\log\frac{|z|}{dt}\biggr]<0.
\end{eqnarray*}
Thus, $\Phi_t$ takes only finitely many positive values and therefore
the maxima
in \textup{(\ref{max})} exist.
\end{pf}

\begin{remark}\label{remark7}
The maximizers $Z_t^{({1})}$,
$Z_t^{({2})}$ and $Z_t^{({3})}$ are in general not uniquely defined.
However, almost surely, if $t_0$ is sufficiently large,
they are uniquely defined, for all $t\in(t_0,\infty)\setminus T$,
where $T$ is a countable random set, and for $t\in T$ it can only happen
that $Z_t^{({1})}=Z_t^{({2})}\neq Z_t^{({3})}$ or $Z_t^{({1})}\neq
Z_t^{({2})}=Z_t^{({3})}$. Thus, the
nonuniqueness only occurs
at the time when the maximal (or the second maximal value) relocates
from one point
to the other. It can be seen from the further
proofs (see Lemma \ref{mon}) that $T$ consists of isolated points.
\end{remark}

To prove Proposition \ref{gap} {below, we need} to analyze the functions
$t\mapsto\Phi_t(Z_t^{({i})})$, $i=1,2,3$, locally. It turns out
that they have
some regularity and that, using rather precise asymptotics for
$|Z_t^{({i})}|$ and
$\Phi_t(Z_t^{({i})})$, we can have good control on their increments.

\begin{lemma}
\label{mon}
Let $\varepsilon>0$. For $i=1,2,3$, almost surely eventually for all
$t$:
\begin{longlist}
\item $\Phi_{t}(Z_{t}^{({i})})> a_t (\log t)^{-\varepsilon
}$ and $\xi(Z_{t}^{({i})})>a_t(\log t)^{-\varepsilon}$;
\item  $t\xi(Z_{t}^{({i})})>|Z_{t}^{({i})}|$;
\item  $r_t(\log t)^{-{1}/{d}-\varepsilon
}<|Z_{t}^{({i})}|< r_t (\log t)^{{1}/({\alpha-d})+\varepsilon}$;
\item  $\Phi_{u}(Z_{u}^{({i})})-\Phi_{t}(Z_{t}^{({i})})\le
\frac{u-t}{t} a_u (\log u)^{{1}/({\alpha-d})+\varepsilon}$ for
all $u>t$;
\item  $u\mapsto\Phi_{u}(Z_{u}^{({i})})$ is increasing on
$(t,\infty)$.
\end{longlist}
\end{lemma}

\begin{pf}As an auxiliary step, let us show that, for any $c>0$ and any
$i\in\mathbb N$,
%
\begin{eqnarray}
\label{mmm}
\xi^{({i-1})}_r
> r^{{d}/{\alpha}}(\log r)^{-c}\qquad \mbox{eventually.}
\end{eqnarray}
Obviously, the distribution of $\xi^{({i-1})}_r$ is given by
\begin{eqnarray*}
\operatorname{Prob}\bigl(\xi^{({i-1})}_r\le x\bigr)=\sum
_{k=0}^{i-1}\pmatrix{\ell_r\cr k} x^{-\alpha k}(1-x^{-\alpha
})^{\ell_r-k},
\end{eqnarray*}
where {$\ell_r\sim\kappa_dr^d$} is the number of points in the ball
$\{|z|\le r\}$,
and $\kappa_d$ is a positive constant.
Using that ${\ell_r\choose k}\le\ell_r^k\sim\kappa_d^kr^{dk}$, we get
\begin{eqnarray*}
&&\operatorname{Prob}\bigl(\xi^{({i-1})}_r\le r^{{d}/{\alpha
}}(\log
r)^{-c}\bigr)\\
&&\qquad\le\bigl(1+o(1)\bigr) \sum_{k=0}^{i-1}\kappa_d^k(\log r)^{c\alpha
k}\bigl(1-r^{-d}(\log r)^{c\alpha}\bigr)^{\ell_r-k}\\
&&\qquad\le\bigl(1+o(1)\bigr) i\kappa_d^{i-1}(\log r)^{c\alpha(i-1)}\bigl(1-r^{-d}(\log
r)^{c\alpha}\bigr)^{\ell_r-i+1}\\
&&\qquad=\exp\bigl\{-\kappa_d(\log r)^{c\alpha}\bigl(1+o(1)\bigr)\bigr\},
\end{eqnarray*}
which is summable along the subsequence $r_n=2^n$. Hence, by the
Borel--Cantelli lemma the inequality \textup{(\ref{mmm})}
holds eventually along $(r_n)_{n\in\mathbb N}$. As $\xi^{({i-1})}_r$ is
increasing, we obtain eventually
\begin{eqnarray*}
\xi^{({i-1})}_r&\ge&\xi^{({i-1})}_{2^{\lfloor\log_2 r\rfloor}}
\ge\bigl(2^{\lfloor\log_2 r\rfloor}\bigr)^{{d}/{\alpha}}\bigl(\log
2^{\lfloor\log_2 r\rfloor}\bigr)^{-c}
\ge2^{-{d}/{\alpha}}r^{{d}/{\alpha}}(\log r-\log2)^{-c}\\
&>& r^{{d}/{\alpha}}(\log r)^{-2c},
\end{eqnarray*}
which is equivalent to \textup{(\ref{mmm})}.

Now we prove parts (i)--(v) of
the lemma. We assume throughout the proof that~$t$ is sufficiently large
to use all statements which hold eventually.

(i) Let $z_1$, $z_2$, $z_3$ be the points where the three largest
values of $\xi$ in
$\{|z|\le r_t(\log t)^{-\varepsilon}\}$ are achieved.
Take $c<\varepsilon(\alpha-d)/(2\alpha)$ and observe that \textup
{(\ref{mmm})} implies for each $i$ eventually
\begin{eqnarray*}
\xi(z_i)>r_t^{{d}/{\alpha}} (\log t)^{-{\varepsilon
d}/{\alpha}}(\log r_t-\varepsilon\log\log t)^{-c}
>a_t (\log t)^{-{\varepsilon d}/{\alpha}-2c}.
\end{eqnarray*}
%
By \cite{HMS08}, Lemma 3.5, we also have
\begin{eqnarray*}
\log\xi(z_i)\le\log\xi^{({1})}_{r_t(\log t)^{-\varepsilon}}<
\log r_t\le(q+1)\log t.
\end{eqnarray*}
We obtain, observing that $t\xi(z_i)>ta_t (\log t)^{-
{\varepsilon d}/{\alpha}-2c}>r_t(\log t)^{-\varepsilon}\ge|z_i|$, that
\begin{eqnarray*}
\Phi_t(z_i)&\ge&\xi(z_i)-\frac{|z_i|}{t}\log\xi(z_i)
>a_t (\log t)^{-{\varepsilon d}/{\alpha}-2c}-\frac
{r_t}{t}(q+1)(\log t)^{1-\varepsilon}\\
&>&a_t(\log t)^{-\varepsilon}
\end{eqnarray*}
as $\frac{\varepsilon d}{\alpha}+2c<\varepsilon$ and $(r_t/t)\log
t=a_t$. Since the inequality is fulfilled
for the three points $z_1, z_2$ and $z_3$, it is also fulfilled for the
maximizers $Z_t^{({1})}$, $Z_t^{({2})}$ and $Z_t^{({3})}$,
completing the proof of the first inequality in (i).
As $\Phi_t(Z_t^{({i})})\neq0$ we must have $\xi(Z_t^{({i})})\ge
|Z_t^{({i})}|/t$, and hence
\begin{eqnarray*}
\xi\bigl(Z_t^{({i})}\bigr)&=& \Phi_{t}\bigl(Z_{t}^{({i})}\bigr)+\frac
{|Z_{t}^{({i})}|}{t}\log\xi\bigl(Z_{t}^{({i})}\bigr)
-\frac{\eta(Z_{t}^{({i})})}{t}\ge\Phi_{t}\bigl(Z_{t}^{({i})}\bigr)+\frac
{|Z_{t}^{({i})}|}{t}\log
\frac{|Z_{t}^{({i})}|}{dt}\\
&>&\Phi_{t}\bigl(Z_{t}^{({i})}\bigr)-d/e.
\end{eqnarray*}
The second inequality in (i) follows now from the lower bound for
$\Phi_{t}(Z_{t}^{({i})})$.\vspace*{-6pt}

\begin{longlist}[(iii)]
\item[(ii)] This is an obvious consequence of (i) as $\Phi
_t(Z_t^{({i})})\neq0$.
\item[(iii)] To prove the upper bound, let us pick $c\in(0,\frac
{\varepsilon(\alpha-d)}{2\alpha})$. Then for each
$z$ such that $|z|\ge r_t(\log t)^{{1}/({\alpha-d})+\varepsilon}$
we obtain by \cite{HMS08}, Lemma 3.5, eventually,
\begin{eqnarray*}
\frac{\xi(z)}{|z|}\le|z|^{{d}/{\alpha}-1}(\log|z|)^{
{1}/{\alpha}+c} \le o(1/t).
\end{eqnarray*}
Hence (ii) implies that $z\neq Z_t^{({i})}$, {which implies the
upper bound on $|Z_t^{({i})}|$.}

To prove the lower bound, suppose that $|Z_t^{({i})}|\le r_t(\log
t)^{-{1}/{d}-\varepsilon}$.
By \cite{HMS08}, Lemma 3.5,
\begin{eqnarray*}
\xi\bigl(Z_t^{({i})}\bigr)\le\bigl|Z_t^{({i})}\bigr|^{{d}/{\alpha}}\bigl(\log
\bigl|Z_t^{({i})}\bigr|\bigr)^{{1}/{\alpha}+c}
\le a_t(\log t)^{-{d\varepsilon}/{\alpha}+2c},
\end{eqnarray*}
which contradicts (i) if we pick $c\in(0,\frac{\varepsilon
d}{2\alpha})$.
\item[(iv)] Let $t$ be large enough so that the previous eventual estimates
hold for all $u\ge t$.
Then, for each $s\in[t,u]$, according to (iii), we have that $\Phi
_{s}(Z_{s}^{({i})})$ is the
$i$th largest value of $\Phi_s$ over {a collection of} finitely many points.
Hence $s\mapsto\Phi_{s}(Z_{s}^{({i})})$ is a continuous piecewise
smooth function.
On the smooth pieces, using again \cite{HMS08}, Lemma 3.5, and (iii)
with $\varepsilon/2$,
we can estimate its derivative by
\begin{eqnarray*}
\frac{{d}}{{ds}} \Phi_{s}\bigl(Z_{s}^{({i})}\bigr)&=&\frac
{|Z_{s}^{({i})}|}{s^2}\log\xi\bigl(Z_{s}^{({i})}\bigr)
-\frac{\eta(Z_{s}^{({i})})}{s^2}
\le\frac{|Z_{s}^{({i})}|}{s^2}\log\bigl|Z_{s}^{({i})}\bigr|^{
{d}/{\alpha}+c}\\
&<&\frac{a_s}{s} (\log s)^{{1}/({\alpha-d})+\varepsilon}.
\end{eqnarray*}
Finally, {we obtain}
\begin{eqnarray*}
\Phi_{u}\bigl(Z_{u}^{({i})}\bigr)-\Phi_{t}\bigl(Z_{t}^{({i})}\bigr)
=\int_t^u\frac{{d}}{{ds}} \Phi_{s}\bigl(Z_{s}^{({i})}\bigr)\,ds
\le\frac{u-t}{t} a_u(\log u)^{{1}/({\alpha-d})+\varepsilon},
\end{eqnarray*}
which completes the proof.
\item[(v)] Using $\eta(z)\le|z| \log d$ in the second, and (i) in
the last step, we see that
\[
\frac{d}{ds} \Phi_{s}\bigl(Z_{s}^{({i})}\bigr) =\frac
{|Z_{s}^{({i})}|}{s^2}\log\xi\bigl(Z_{s}^{({i})}\bigr)
-\frac{\eta(Z_{s}^{({i})})}{s^2} \ge\frac{|Z_{s}^{({i})}|}{s^2}\log
\frac{\xi(Z_{s}^{({i})})}{d}>0,
\]
eventually for all $t$.\quad\qed
\end{longlist}
\noqed\end{pf}

\subsection{Lower bound for $\Phi_t(Z_t^{({1})})-\Phi
_t(Z_t^{({3})})$}\label{sec-Phi1-Phi3}

In this section we prove that $\Phi_t(Z_t^{({1})})$ and $\Phi
_t(Z_t^{({3})})$ are well separated from each other.
The crucial estimate for this is provided in Lemma \ref{lem-Phi1-Phi3}.

First, it is important to make the density of the random variable
$\Phi_t(z)$ explicit. Observe that, on the set
$\{t\xi(z)\geq z\}$, the event $\{\Phi_t(z)<x\}$ has the form $\{\chi
_a(\xi(z))\leq x-\eta(z)/t\}$, where we
abbreviated $a=|z|/t$ and introduced the map $\chi_a(x)=x-a\log x$.
Note that $\chi_a$ is an increasing bijection
from $[a,\infty)$ to $[a-a\log a,\infty)$, hence on $\{t\xi(z)\geq
z\}$ we can describe $\{\Phi_t(z)<x\}$ using
the inverse function $\psi_a\dvtx[a-a\log a,\infty)\to[a,\infty)$
of $\chi_a$. In order to also include the
complement of $\{t\xi(z)\geq z\}$, we extend $\psi_a$ to a function
$\mathbb R\to[a,\infty)$ by putting $\psi_a(x)=a$ for $x<a-a\log a$.
Then we have, for each $t$, $z$ and $x>0$,
%
\begin{eqnarray}
\label{ppsi}
\{\Phi_t(z)\le x\}
=\bigl\{\xi(z)\le\psi_{{|z|}/{t}}\bigl(x-\eta(z)/t\bigr)\bigr\}.
\end{eqnarray}

\begin{lemma}\label{lem-Phi1-Phi3}
\label{prob}
Fix $\beta>1+\frac{1}{\alpha-d}$ and let $\lambda_t=(\log
t)^{-\beta}$.
Then there exists a~constant $c>0$ such that
\begin{eqnarray*}
\operatorname{Prob}\bigl(\Phi_t\bigl(Z_t^{({1})}\bigr)-\Phi_t\bigl(Z_t^{({3})}\bigr)\le
2a_t\lambda_t\bigr)\le c\lambda_t^2\qquad
\mbox{for } t>0.
\end{eqnarray*}
\end{lemma}

\begin{pf}
This proof, though tedious, is fairly standard and is carried out in
\textit{four steps}. In the first step, we
show that there exists a constant $C_1>0$ such that, for all
sufficiently large $t$, and all $s\ge(\log t)^{-1/2}$,
%
\begin{eqnarray}\label{ersterschritt}\qquad
&&\operatorname{Prob} \bigl( \Phi_t\bigl(Z_t^{({1})}\bigr) \in d(a_t s), \Phi
_t\bigl(Z_t^{({1})}\bigr) -
\Phi_t\bigl(Z_t^{({3})}\bigr)
\le2a_t \lambda_t \bigr)\nonumber\\[-8pt]\\[-8pt]
&&\qquad \le C_1 a_t^3 \lambda_t^2 \operatorname{Prob}\bigl( \Phi
_t\bigl(Z_t^{({1})}\bigr) \le a_ts\bigr)
\Biggl[ \sum_{z\in\mathbb Z^d} \biggl(a_ts + \frac{|z|}{t} \log\frac
{|z|}{dt}\biggr)^{-\alpha-1} \Biggr]^3 \,ds .\nonumber
\end{eqnarray}
In the second step we evaluate the infinite sum and show that there
exists $C_2>0$ such that
%
\begin{equation}\label{zweiterschritt}
\sum_{z\in\mathbb Z^d} \biggl(a_ts + \frac{|z|}{t} \log\frac
{|z|}{dt}\biggr)^{-\alpha-1}
\le C_2 a_t^{-1} s^{d-\alpha-1} .
\end{equation}
To bound the right-hand side of \textup{(\ref{ersterschritt})}
further, we show
in the third
step that there exists a constant $C_3>0$ such that, for all $(\log
t)^{-1/2}\le s\le1$,
%
\begin{equation}\label{dritterschritt}
\operatorname{Prob}\bigl( \Phi_t\bigl(Z_t^{({1})}\bigr) \le a_ts \bigr) \le
\exp\{ -
C_3 s^{d-\alpha} \} .
\end{equation}
In the fourth step we combine these three equations and integrate
over $s$ to get the result.

For the \textit{first step} we use independence to obtain
%
\begin{eqnarray}\label{bigdisplay}
&&\operatorname{Prob}\bigl( \Phi_t\bigl(Z_t^{({1})}\bigr) \in d(a_t s), \Phi
_t\bigl(Z_t^{({1})}\bigr) -
\Phi_t\bigl(Z_t^{({3})}\bigr)
\le2a_t \lambda_t \bigr)\nonumber\\[-1pt]
&&\qquad \le\mathop{\mathop{\sum}_{z_1, z_2, z_3\in\mathbb Z^d}}_{\mathrm
{distinct}}
\operatorname{Prob}\bigl( \Phi_t(z_1) \in d(a_t s);\nonumber\\[-4pt]
&&\hspace*{103.1pt} \Phi_t(z_i)
\in
a_t[s-2\lambda_t,s] \mbox{ for }i=2,3; \nonumber\\[-1pt]
&&\hspace*{103.1pt}\hspace*{20.2pt} \Phi_t(z) \le a_t s \mbox{ for }
z\notin\{z_1,z_2,z_3\} \bigr)\\[-2pt]
&&\qquad \le\Biggl( \sum_{z\in\mathbb Z^d} \frac{\operatorname{Prob}( \Phi
_t(z) \in d(a_t
s))}{\operatorname{Prob}(\Phi_t(z) \le a_t s)}\Biggr)\nonumber\\[-1pt]
&&\qquad\quad{}\times
\Biggl(\sum_{z\in\mathbb Z^d} \frac{\operatorname{Prob}( \Phi_t(z)
\in a_t[s-2\lambda
_t,s])}{\operatorname{Prob}(\Phi_t(z) \le a_t s)} \Biggr)^2 \nonumber\\
&&\qquad\quad{} \times\prod_{z\in\mathbb Z^d} \operatorname
{Prob}\bigl( \Phi
_t(z) \le a_t s\bigr) .\nonumber
\end{eqnarray}
All the denominators in \textup{(\ref{bigdisplay})} converge to one,
uniformly in $z$
and $s\ge(\log t)^{-1/2}-2\lambda_t$.
Indeed, by \textup{(\ref{ppsi})}, we get
\begin{eqnarray*}
\operatorname{Prob}\bigl(\Phi_t(z) \le a_t s \bigr) &=&
\operatorname{Prob}\biggl( \xi(z) \le\psi_{{|z|}/{t}}\biggl(a_t
s-\frac{\eta(z)}{t} \biggr)\biggr) \\[-1pt]
&\ge& \operatorname{Prob}\biggl( \xi(z) \le a_ts - \frac{|z|}{t} \log
d + \frac{|z|}{t} \log\frac{|z|}{t} \biggr)\\[-1pt]
&\ge&\operatorname{Prob}\biggl( \xi(z) \le a_ts - \frac{d}{e} \biggr) \ge1 + o(1),
\end{eqnarray*}
using that $\eta(z) \le|z| \log d$, $x \log(x/d) \ge-d/
e$, and
$\psi_a(x)\ge x + a \log a$ (with $a=|z|/t$), where the latter is
obvious for $x\le a-a\log a$
and follows from $\psi_a(x)=x+a\log\psi_a(x)\ge x+a\log a$ otherwise.

Further, we use \textup{(\ref{ppsi})} to observe that, by a coordinate
transformation, the density of $\Phi_t(z)$ at $x$ is given as
\[
\psi'_{{|z|}/{t}}\biggl(x - \frac{\eta(z)}{t}\biggr) \alpha
\biggl( \psi_{{|z|}/{t}}\biggl( x - \frac{\eta(z)}{t} \biggr)
\biggr)^{-\alpha-1},\qquad
\mbox{if } x - \frac{\eta(z)}{t}> \frac{|z|}{t}- \frac
{|z|}{t} \log\frac{|z|}{t}.
\]
If $t$ is large enough, the latter condition is satisfied for $x=a_ts$,
all $z$
and $s\ge(\log t)^{-1/2}-2\lambda_t$, and moreover, using again $\psi
_a(x)\ge x + a \log a$, we have
\[
\psi_{{|z|}/{t}}\biggl(a_ts - \frac{\eta(z)}{t}\biggr)
\ge a_t s - \frac{\eta(z)}{t} + \frac{|z|}{t} \log\frac{|z|}{t}
\ge a_ts
+ \frac{|z|}{t} \log\frac{|z|}{dt} .
\]
Hence, if $t$ is big enough to satisfy $a_t[(\log t)^{-1/2}-2\lambda
_t]>t^{q/2}$ we get
\[
\frac{t}{|z|} \psi_{{|z|}/{t}}\biggl(a_ts - \frac{\eta(z)}{t}\biggr)
\ge
\frac{1}{|z|} t^{1+q/2} + \log\frac{|z|}{dt} \ge\min_{r>0}
\{ t^{q/2}r - \log(rd) \}
= \log\frac{et^{q/2}}{d}.
\]
Differentiating the equality $\psi_a(x)-a\log\psi_a(x)=x$ with
respect to $x$, for $x>a$, we obtain
$\psi'_a(x)=(1-a/\psi_a(x))^{-1}$.
This implies that, as $t\uparrow\infty$,
\begin{eqnarray}
\psi'_{{|z|}/{t}}\biggl(a_ts - \frac{\eta(z)}{t}\biggr) = \biggl(1-
\frac{|z|}{t}\Big/\psi_{{|z|}/{t}}\biggl(a_ts - \frac{\eta
(z)}{t}\biggr) \biggr)^{-1}
\longrightarrow1\nonumber\\
\eqntext{\mbox{uniformly in $z$ and $s$.}}
\end{eqnarray}
Hence
%
\begin{equation}\label{densityformula}
\operatorname{Prob}\bigl( \Phi_t(z) \in d(a_t s)\bigr) \le\bigl(\alpha
+o(1)\bigr) a_t
\biggl( a_ts
+ \frac{|z|}{t} \log\frac{|z|}{dt} \biggr)^{-\alpha-1}\,ds .
\end{equation}
Integrating \textup{(\ref{densityformula})} over the interval
$[s-2\lambda
_t,s]$ yields
\[
\operatorname{Prob}\bigl( \Phi_t(z) \in a_t [s-2\lambda_t,s]\bigr)
\le\bigl(\alpha+o(1)\bigr) 2a_t\lambda_t \biggl( a_t(s-2\lambda_t)
+ \frac{|z|}{t} \log\frac{|z|}{dt} \biggr)^{-\alpha-1}.
\]
Using that $x \log(x/d) \ge-d/e$ we obtain
\begin{eqnarray*}
\frac{a_t(s-2\lambda_t) + {|z|}/{t} \log({|z|}/{dt})}
{a_ts + ({|z|}/{t}) \log({|z|}/{dt})} &=&
1- \frac{2a_t\lambda_t}{a_ts + ({|z|}/{t}) \log({|z|}/{dt})}\\
&\ge&1 - \frac{2a_t\lambda_t}{a_ts -d/e}\ge1+o(1),
\end{eqnarray*}
%
hence, uniformly in $z\in\mathbb Z^d$ and $s\ge(\log
t)^{-1/2}-2\lambda_t$,
%
\begin{equation}\label{alsotrue}
\operatorname{Prob}\bigl( \Phi_t(z) \in a_t [s-2\lambda_t,s]\bigr)
\le\bigl(\alpha+o(1)\bigr) 2a_t\lambda_t \biggl( a_ts + \frac{|z|}{t}
\log\frac{|z|}{dt} \biggr)^{-\alpha-1}.\hspace*{-38pt}
\end{equation}
Inserting \textup{(\ref{densityformula})} and \textup{(\ref
{alsotrue})} in \textup{(\ref{bigdisplay})} and estimating
all denominators uniformly by a~constant factor yields \textup{(\ref
{ersterschritt})}.

In the \textit{second step} we estimate the infinite sum in \textup{(\ref
{zweiterschritt})}.
Recalling that $r_t^d=a_t^\alpha$ and that the number of points in the
{ball $\{|z|\le r\}$}
is equal to $\kappa_d r^d (1+o(1))$ we obtain
%
\begin{eqnarray}\label{smallpart}\quad
\sum_{|z|\le r_t/\log t} \biggl(a_ts + \frac{|z|}{t} \log\frac
{|z|}{dt}\biggr)^{-\alpha-1}
& \le&\sum_{|z|\le r_t/\log t} \biggl(a_ts -\frac{d}{e}
\biggr)^{-\alpha-1}\nonumber\\
& \le&\bigl(\kappa_d+o(1)\bigr) \frac{r_t^d}{(a_ts)^{\alpha+1}[\log t]^d}\\
& =&
o( a_t^{-1} s^{d-\alpha-1}).\nonumber
\end{eqnarray}
We have $\log\frac{|z|}{dt} \ge(1+o(1)) q [\log t]$ uniformly
over all $z\in\mathbb Z^d$
with $|z|\ge r_t/\log t$. Therefore,
\begin{eqnarray*}
&&
\sum_{|z|\ge r_t/\log t} \biggl(a_ts + \frac{|z|}{t} \log
\frac{|z|}{dt}\biggr)^{-\alpha-1}\\
&&\qquad \le\bigl(1+o(1)\bigr) (a_t s)^{-\alpha-1} \sum_{|z|\ge r_t/\log t}
\biggl(1 + q
\frac{|z|}{r_t s} \log t\biggr)^{-\alpha-1} \\
&&\qquad =\bigl(1+o(1)\bigr) \frac{(r_t s)^d}{(a_t s)^{\alpha+1}}
\int_{\mathbb R^d} (1 + q |x|)^{-\alpha-1} \,dx
\le C_2a_t^{-1}s^{d-\alpha-1}.
\end{eqnarray*}
%
Combining this with \textup{(\ref{smallpart})} yields \textup{(\ref
{zweiterschritt})}.

In the \textit{third step} we show \textup{(\ref{dritterschritt})} by a direct
calculation.
First, let us show that for $\varepsilon>0$
%
\begin{eqnarray}
\label{kkk}
\psi_{{|z|}/{t}}(a_ts)\le a_ts+\frac{|z|}{t} (q+\varepsilon
)\log t
\end{eqnarray}
for all large $t$, $\frac{r_t}{\log t}\le|z|\le r_t\log t$ and $(\log
t)^{-{1}/{2}}\le s\le1$.
{By definition,
\[
\psi_{{|z|}/{t}}(a_ts)= a_t s + \frac{|z|}{t} \log\psi_{
{|z|}/{t}}(a_t s),
\]
hence
it suffices to prove that $\psi_{{|z|}/{t}}(a_ts)\le
t^{q+\varepsilon}$.}

Assume this is false for some large $t$, $z$ and $s$.
Then using the monotonicity of $x\mapsto x-a\log x$ for $x\ge a$, we obtain
\begin{eqnarray*}
t^{q+{\varepsilon}/{2}}&\ge& a_t\ge a_ts=\psi_{
{|z|}/{t}}(a_ts)-\frac{|z|}{t}\log\psi_{{|z|}/{t}}(a_ts)
\ge t^{q+\varepsilon}-\frac{|z|}{t} (q+\varepsilon)\log t\\
&\ge&
t^{q+\varepsilon}-t^{q+{\varepsilon}/{2}},
\end{eqnarray*}
which is a contradiction.
Now we can compute
\begin{eqnarray*}
\operatorname{Prob}\bigl( \Phi_t\bigl(Z_t^{({1})}\bigr) \le a_ts \bigr)
&=& \prod_{z\in\mathbb Z^d} \operatorname{Prob}\bigl( \Phi_t(z) \le
a_ts \bigr)\\
&\le& \prod_{{r_t}/({\log t})\le|z|\le r_t\log t}
\operatorname{Prob}\biggl(\xi(z) \le\psi_{{|z|}/t}\biggl(a_t s-\frac
{\eta
(z)}{t}\biggr)\biggr)\\
&\le&\prod_{{r_t}/({\log t})\le|z|\le r_t\log t}
\operatorname{Prob}\biggl(\xi(z) \le a_ts+\frac{|z|}{t}
(q+\varepsilon)\log
t\biggr)
\end{eqnarray*}
using \textup{(\ref{ppsi})}, $\psi_{{|z|}/t}(a_t s-\frac{\eta(z)}{t})
\le
\psi_{{|z|}/t}(a_t s)$ and \textup{(\ref{kkk})}.
Inserting the explicit form of the distribution
function we get
\begin{eqnarray*}
&&\operatorname{Prob}\bigl( \Phi_t\bigl(Z_t^{({1})}\bigr) \le a_ts
\bigr)\\
&&\qquad\le
\exp\Biggl\{ -\bigl(1+o(1)\bigr) \sum_{{r_t}/({\log t})\le
|z|\le r_t\log t}
\biggl( a_t s + \frac{|z|}{t} (q+\varepsilon)\log t\biggr)^{-\alpha
}\Biggr\}\\
&&\qquad\le\exp\biggl\{ -\bigl(1+o(1)\bigr) s^{d-\alpha}\int_{\mathbb R^d} \frac{d
u}{(1+(q+\varepsilon)|u|)^\alpha}\biggr\}
\end{eqnarray*}
using a Riemann sum approximation as in the second step. This proves
\textup{(\ref{dritterschritt})}.

Coming to the \textit{fourth step}, we now use \textup{(\ref{ersterschritt})},
\textup{(\ref{zweiterschritt})} and \textup{(\ref{dritterschritt})}
to get
\begin{eqnarray*}
&&\operatorname{Prob}\bigl( \Phi_t\bigl(Z_t^{({1})}\bigr) - \Phi
_t\bigl(Z_t^{({3})}\bigr) \le2a_t
\lambda_t \bigr)\\
&&\qquad \le\operatorname{Prob}\bigl( \Phi_t\bigl(Z_t^{({1})}\bigr) \le a_t(\log
t)^{-1/2} \bigr)
\\
&&\qquad\quad{} + \int_{(\log t)^{-1/2}}^\infty\operatorname{Prob}
\bigl(\Phi
_t\bigl(Z_t^{({1})}\bigr)
\in d(a_t s), \Phi_t\bigl(Z_t^{({1})}\bigr) - \Phi_t\bigl(Z_t^{({3})}\bigr) \le
2a_t \lambda_t \bigr)\\
&&\qquad \le\exp\bigl\{ - C_3 (\log t)^{({\alpha-d})/{2}} \bigr\}\\
&&\qquad\quad{} + C_1C_2^3 \lambda_t^2 \biggl[ \int_{(\log t)^{-1/2}}^1
\frac{\exp\{ - C_3 s^{d-\alpha} \}\,ds}{s^{3(\alpha-d+1)}} +
\int_1^\infty\frac{ds}{s^{3(\alpha-d+1)}} \biggr] .
\end{eqnarray*}
The first term is $o(\lambda_t^2)$ by choice of $\lambda_t$, and the
expression in the square bracket
is bounded by an absolute constant. This completes the proof.
\end{pf}

Now we turn the estimate of Lemma \ref{lem-Phi1-Phi3} into an almost
sure bound.

\begin{prop}
\label{gap}
Almost surely, eventually for all $t$,
\begin{eqnarray*}
\Phi_t\bigl(Z_t^{({1})}\bigr)-\Phi_t\bigl(Z_t^{({3})}\bigr) \ge a_t\lambda_t.
\end{eqnarray*}
\end{prop}

\begin{pf}Let $\varepsilon\in(0,2\beta-1)$ be such that $\beta
>1+\frac{1}{\alpha-d}+2\varepsilon$
and let $t_n=e^{n^{\gamma}}$, where $\gamma\in(\frac
{1}{2\beta},\frac{1}{2\beta-\varepsilon})$. Note that $\gamma<1$.
Since $\lambda_{t_n}^2=n^{-2\gamma\beta}$ is summable, Lemma~\ref
{prob} and the Borel--Cantelli
lemma imply that
\begin{eqnarray*}
\Phi_{t_n}\bigl(Z_{t_n}^{({1})}\bigr)-\Phi_{t_n}\bigl(Z_{t_n}^{{({3})}}\bigr)\ge
2a_{t_n}\lambda_{t_n}\qquad \mbox{eventually for all }n.
\end{eqnarray*}
For each $t\in[t_n,t_{n+1})$ we obtain by Lemma \ref{mon}(iv,v)
%
\begin{eqnarray}\label{ac}
&&\Phi_t\bigl(Z_t^{({1})}\bigr)-\Phi_t\bigl(Z_t^{({3})}\bigr)\nonumber\\
&&\qquad\ge\Phi_{t_n}\bigl(Z_{t_n}^{({1})}\bigr)-\Phi
_{t_{n+1}}\bigl(Z_{t_{n+1}}^{({3})}\bigr)\nonumber\\[-8pt]\\[-8pt]
&&\qquad=\bigl[\Phi_{t_n}\bigl(Z_{t_n}^{({1})}\bigr)-\Phi_{t_n}\bigl(Z_{t_n}^{({3})}\bigr)\bigr]
-\bigl[\Phi_{t_{n+1}}\bigl(Z_{t_{n+1}}^{({3})}\bigr)-\Phi
_{t_n}\bigl(Z_{t_n}^{({3})}\bigr)\bigr]\nonumber\\
&&\qquad\ge2a_{t_n}\lambda_{t_n}-\frac{t_{n+1}-t_n}{t_n} a_{t_{n+1}}(\log
t_{n+1})^{{1}/({\alpha-d})+\varepsilon}.\nonumber
\end{eqnarray}
Notice that eventually
\begin{eqnarray*}
\frac{t_{n+1}-t_n}{t_n}&=&e^{(n+1)^{\gamma}-n^{\gamma
}}-1=\gamma
n^{\gamma-1}\bigl(1+o(1)\bigr)
=\gamma(\log t_n)^{({\gamma-1})/{\gamma}}\bigl(1+o(1)\bigr)\\
&\le&(\log
t_n)^{-2\beta+1+\varepsilon}.
\end{eqnarray*}
Denote by $n(t)$ the integer such that $t\in[t_{n(t)},t_{n(t)+1})$.
Since $t_{n+1}/t_n\to1$ we have
$t_{n(t)}\sim t$ and $t_{n(t)+1}\sim t$. Substituting this and the last
estimate into \textup{(\ref{ac})},
we obtain
\begin{eqnarray*}
\Phi_t\bigl(Z_t^{({1})}\bigr)-\Phi_t\bigl(Z_t^{({3})}\bigr)
&\ge&2a_t\lambda_t\bigl(1+o(1)\bigr)-a_t(\log t)^{{1}/({\alpha-d})-2\beta
+1+2\varepsilon}\bigl(1+o(1)\bigr)\\
&\ge& a_t\lambda_t
\end{eqnarray*}
eventually since $\lambda_t=(\log t)^{-\beta}$ and
$(\log t)^{{1}/({\alpha-d})-\beta+1+2\varepsilon}=o(1)$, which makes
the second term negligible {compared} to the first one.
\end{pf}

\section{Total mass of the solution}\label{s_mass}

In this section we show that the total mass~$U(t)$ of the solution can
be well approximated by
$\exp\{t \Phi_t(Z_t^{({1})})\}$. The main tool is the Feynman--Kac formula
in \textup{(\ref{FKform})} and a technical lemma
provided in Section \ref{sec-Jtesti}. In Section \ref{sec-LB} we
prove the lower bound for
$\frac1t\log U(t)$. In Section \ref{sec-U12345} we split the set of
all paths into five path classes,
four of which turn out to give negligible contribution to the
Feynman--Kac formula for $U(t)$. In Section \ref{sec-UB} we show that
the remaining class yields a useful upper bound for $\frac1t\log U(t)$.

\subsection{A technical lemma}\label{sec-Jtesti}

We bound contributions to the Feynman--Kac formula for $U(t)$ coming
from path classes that
are defined according to their number of steps and the maximum along
their path. Denote by $J_t$ the number of jumps
of the random walk $(X_s\dvtx s\ge0)$ {up to time $t$.} Recall the
notation from the beginning of Section \ref{s_pre}
and let $H=(H_t)_{t>0}$ be some family of sets $H_t\subset F_t$, and
let $h=(h_t)_{t>0}$ be some family of functions
$h_t\dvtx\mathbb Z^d\to\mathbb N_0$. Denote by
\begin{eqnarray*}
&&U_{H,h}(t)=\mathbb E_0\biggl[\exp\biggl\{\int_0^t\xi(X_s)\,ds\biggr\}\\
&&\hspace*{63pt}{}\times
{\mathbh1}\biggl\{\exists z\in F_t\setminus H_t\dvtx
\max_{s\in[0,t]}\xi(X_s)=\xi(z), { h_t(z) \le J_t \le
r_tg_t} \biggr\}\biggr]
\end{eqnarray*}
the contribution to the total mass that comes from paths which attain
their maximal potential value
in some $z\in F_t\setminus H_t$ with step number in $\{h_t(z),\ldots,
{\lfloor r_tg_t \rfloor}\}$.

\begin{lemma}
\label{l10}
There is $\delta>0$ such that, almost surely, for $t\to\infty$:
\begin{longlist}[{(a)}]
\item[{(a)}]
$\frac{1}{t}\log U_{H,h}(t)\le\max_{z\in F_t\setminus H_t}\{
\Phi_t(z)
+\frac{1}{t}\max_{n\ge h_t(z)}[\eta(n,z)-\eta(z)-\frac
{n-|z|}{2}\times\break\log\xi(z)]\}+O(t^{q-\delta})$,
\item[{(b)}]
$\frac{1}{t}\log U_{H,h}(t)\le\max_{z\in F_t\setminus H_t}\Phi
_t(z)+O(t^{q-\delta}).$
\end{longlist}
\end{lemma}

\begin{pf} {We write} $U_{H,h}(t)$ as
%
\begin{eqnarray}
\label{tz}
U_{H,h}(t)=\sum_{z\in F_t\setminus H_t}U_{H,h}(t,z),
\end{eqnarray}
where we define, {for any $z\in\mathbb Z^d$,}
\begin{eqnarray*}
U_{H,h}(t,z)=\mathbb E_0\biggl[\exp\biggl\{\int_0^t\xi(X_s)\,ds\biggr\}
{\mathbh1}\biggl\{\max_{s\in[0,t]}\xi(X_s)=\xi(z), h_t(z)
\le
J_t\le r_tg_t\biggr\}\biggr].
\end{eqnarray*}
Denote by
\begin{eqnarray*}
&&{\mathcal P}_{n,z}=\biggl\{y=(y_0,y_1,\ldots,y_n)\in(\mathbb Z^d)^{n+1}
\dvtx y_0=0, |y_i-y_{i-1}|=1,\\
&&\hspace*{203pt}
\max_{0\le i\le n}\xi(y_i)=\xi(z)\biggr\}
\end{eqnarray*}
the set of all discrete time paths in $\mathbb Z^d$ of length $n$
starting at
the origin,
going through $z$, such that the maximum of the potential over the path
is attained at $z$.
Let $(\tau_i)_{i\in\mathbb N_0}$
be a sequence of independent exponentially distributed random variables
with parameter
$2d$. Denote by $\mathsf{E}$ the expectation with respect to $(\tau_i)$.
Averaging over all random paths following the same path $y$
(with individual timings) we obtain
%
\begin{eqnarray}
\label{tzy}
U_{H,h}(t,z)=\sum_{n=h_t(z)}^{{\lfloor}r_tg_t{\rfloor}}\sum_{y\in
{\mathcal P}_{n,z}}U_{H,h}(t,z,y),
\end{eqnarray}
where
\begin{eqnarray*}
&&U_{H,h}(t,z,y)\\
&&\qquad=(2d)^{-n} \mathsf{E}
\Biggl[\exp\Biggl\{\sum_{i=0}^{n-1}\tau_i\xi(y_{i})+
\Biggl(t-\sum_{i=0}^{n-1}\tau_i\Biggr)\xi(y_n)\Biggr\}{\mathbh1}\Biggl\{
\sum_{i=0}^{n-1}\tau_i\leq t<\sum_{i=0}^{n}\tau_i\Biggr\}\Biggr].
\end{eqnarray*}
Note that, as $y$ can have self-intersections,
some of the values of $\xi$ over $y$ may coincide. We would like to
avoid the
situation when the maximum of $\xi$ over $y$ is taken at more than one point.
Therefore, for each path $y$, we slightly change the
potential {over $y$.} Namely, we denote by
$i(y)=\min\{i \dvtx\xi(y_i)=\xi(z)\}$
the index of the first point where the maximum of the potential over
the path is attained. Then we define the modified version of the
potential $\xi^{y} \dvtx\{0,\ldots,n\}\to\mathbb R$ by
\begin{eqnarray*}
\xi^y_i=\cases{
\xi(y_i),&\quad if $i\neq i(y)$,\cr
\xi(y_i)+1,&\quad if $i=i(y)$.}
\end{eqnarray*}
Repeating the computations (4.16) and (4.17) from \cite{HMS08} we obtain
%
\begin{eqnarray}
\label{product}\qquad
U_{H,h}(t,z,y)&\le& e^{t\xi^y_{i(y)}-2dt}\prod_{i\neq
i(y)}\frac
{1}{\xi^y_{i(y)}-\xi^y_i}
\le e^{t\xi(z)}\prod_{i=1}^n\frac{1}{1+\xi(z)-\xi(y_i)}.
\end{eqnarray}
Let us now find a lower bound for the number of {sites on the path
where} the potential is
small compared to $\xi(z)$ { or, more precisely, we estimate the
number of indices $1\le i \le n$}
such that $\xi(y_i)\in G_t^{\mathrm{c}}$. First, we erase loops that
the path $y$ may have made
before reaching $z$ for the first time and extract from $(y_0,\ldots
,y_{i(y)})$ a self-avoiding
path $(y_{i_0},\ldots,y_{i_{l(y)}})$ starting at the origin, ending at
$z$ and having length $l(y)\ge|z|$,
where we take $i_0=0$ and
\begin{eqnarray*}
i_{j+1}=\min\{i\dvtx y_l\neq y_{i_j}\ \forall l\in[i,i(y)]\}.
\end{eqnarray*}
Since this path visits $l(y)$ different points, at least $l(y)-m_t$ of
them belong to $G_t^{\mathrm{c}}$.
By Lemma \ref{asymp}(ii) we have $|z|>t^{(q+1)\sigma+c}>m_t$ and
hence $l(y)-m_t$
is eventually positive. Second, for each $0\le j\le l(y)-1$, consider
the path $(y_{i_j+1},\ldots,y_{i_{j+1}-1})$, which was
removed during erasing the $j$th loop. Obviously, it contains {an even
number $i_{j+1}-{i_j}-1$}
of steps,\vspace*{1pt} as $y_{i_j}=y_{i_{j+1}-1}$ and $y_{i_j}$ and $y_{i_{j+1}-1}$
are neighbors.
Notice that, as $G_t$ is totally disconnected by Lemma \ref{asymp}(iii),
at least half of the steps,
$(i_{j+1}-{i_j}-1)/2$, belong to $G_t^{\mathrm{c}}$.
Third,
consider the remaining piece $(y_{i(y)+1},\ldots,y_n)$. Again, since
$G_t$ is totally disconnected,
there will be at least $(n-i(y))/2$ points belonging to $G_t^{\mathrm{c}}$.
Summing up these three observations, we obtain that $y$ makes at least
\begin{eqnarray*}
l(y)-m_t+\sum_{j=0}^{l(y)-1}\frac{{i_{j+1}}-{i_j}-1}{2}+\frac{n-i(y)}{2}
&=&l(y)-m_t+\frac{n-l(y)}{2}\\
&\ge&|z|-m_t+\frac{n-|z|}{2}
\end{eqnarray*}
steps that belong to $G_t^{\mathrm{c}}$.

Now we can continue estimating $U_{H,h}(t,z,y)$. Recall that the
potential is at most $\xi^{({m_t})}_{r_tg_t}$
on the set $G_t^{\mathrm{c}}$. If we drop the terms corresponding to
the points from $G_t$
in \textup{(\ref{product})}, we obtain
\begin{eqnarray*}
U_{H,h}(t,z,y) 
\le e^{t\xi(z)}\bigl[\xi(z)-\xi^{({m_t})}_{r_tg_t}
\bigr]^{-(|z|-m_t+({n-|z|})/{2})}.
\end{eqnarray*}
Substituting this into \textup{(\ref{tzy})} and using $|{\mathcal
P}_{n,z}|\le
N(n,z)$, we get
\begin{eqnarray*}
&&\hspace*{-4pt}\frac{1}{t}\log U_{H,h}(t,z)\\
&&\hspace*{-4pt}\qquad\le\frac{1}{t}\log\sum_{n=h_t(z)}^{r_tg_t}\sum_{y\in{\mathcal P}_{n,z}}
e^{t\xi(z)}\bigl[\xi(z)-\xi^{({m_t})}_{r_tg_t}
\bigr]^{{-(|z|-m_t+({n-|z|})/{2})}}\\
&&\hspace*{-4pt}\qquad\le\frac{1}{t}\log\max_{h_t(z)\vee|z|\le n\le r_tg_t}\bigl\{
N(n,z)
e^{t\xi(z)}\bigl[\xi(z)-\xi^{({m_t})}_{r_tg_t}
\bigr]^{{-(|z|-m_t+({n-|z|})/{2})}}\bigr\}+o(1)\\
&&\hspace*{-4pt}\qquad=\max_{h_t(z)\vee|z|\le n\le r_tg_t}\biggl\{\xi(z)+\frac{\eta
(n,z)}{t}\\
&&\hspace*{-4pt}\hspace*{100pt}{}-\frac{1}{t}
\biggl[|z|-m_t +\frac{n-|z|}{2}\biggr]\log\bigl(\xi(z)-\xi
^{({m_t})}_{r_tg_t}\bigr)\biggr\}+o(1).
\end{eqnarray*}
In order to simplify the expression under the maximum, we decompose
\begin{eqnarray*}
&&\biggl[|z|-m_t + \frac{n-|z|}{2}\biggr]\log\bigl(\xi(z)-\xi
^{({m_t})}_{r_tg_t}\bigr)\\
&&\qquad =\biggl[|z|+\frac{n-|z|}{2}\biggr]\log\xi(z)
+\biggl[|z|-m_t+\frac{n-|z|}{2}\biggr]\log\biggl(1-\frac{\xi
^{({m_t})}_{r_tg_t}}{\xi(z)}\biggr)\\
&&\qquad\quad{} -m_t\log\xi(z)
\end{eqnarray*}
and {show} that the last two terms are negligible.
Indeed, for the second term, we use Lemma \ref{asymp}(i)
{in the second step} to obtain, for
each $\delta<c$,
\begin{eqnarray*}
\biggl|\biggl[|z|-m_t+\frac{n-|z|}{2}\biggr]\log\biggl(1-\frac{\xi
^{({m_t})}_{r_tg_t}}{\xi(z)}\biggr)\biggr|
&\le& n\biggl|\log\biggl(1-\frac{\xi^{({m_t})}_{r_tg_t}}{\xi
^{({k_t})}_{r_tg_t}}\biggr)\biggr|
\le r_tg_t t^{-c}\\
&=&O(t^{q+1-\delta})
\end{eqnarray*}
uniformly for all $n\ge|z|$. For the third term, we use \cite{HMS08}, Lemma
3.5, and obtain
$\log\xi(z)\le O(\log t)$ uniformly for all $|z|\le r_tg_t$. For
$\delta<(q+1)(1-\sigma)$ this implies that
\begin{eqnarray*}
m_t\log\xi(z)\le O((r_tg_t)^{\sigma}\log t)=O(t^{q+1-\delta}).
\end{eqnarray*}
Hence, there is a small positive $\delta$ such that
%
\begin{eqnarray}\label{nadoelo}\hspace*{26pt}
&&\frac{1}{t}\log U_{H,h}(t,z)\nonumber\\
&&\qquad\le
\max_{h_t(z)\vee|z|\le n\le r_tg_t}\biggl\{\xi(z)+\frac{\eta
(n,z)}{t}-\frac{1}{t}
\biggl[|z|+\frac{n-|z|}{2}\biggr]\log\xi(z)\biggr\}\nonumber\\
&&\qquad\quad{} + O(t^{q-\delta})\\
&&\qquad=\biggl[\xi(z)+\frac{\eta(z)}{t}-\frac{|z|}{t}\log\xi(z)\biggr]\nonumber\\
&&\qquad\quad{} +\frac{1}{t}\max_{h_t(z)\vee|z|\le n\le r_tg_t}
\biggl\{\eta(n,z)-\eta(z)-
\frac{n-|z|}{2}\log\xi(z)\biggr\}+O(t^{q-\delta}).\nonumber\!\!\!\!\!
%
\end{eqnarray}
To prove (a), observe that for each $z\in F_t$ we have $\xi
(z)>ed$.
Hence either $t\xi(z)\ge|z|$ or 
\begin{eqnarray*}
\xi(z)+\frac{\eta(z)}{t}-\frac{|z|}{t}\log\xi(z)
\le\xi(z)-\frac{|z|}{t}\log\frac{\xi(z)}{d}\le\xi(z)
\biggl[1-\log\frac{\xi(z)}{d}\biggr]<0.
\end{eqnarray*}
In any case we obtain, using \textup{(\ref{tz})} and \textup{(\ref
{nadoelo})},
\begin{eqnarray*}
\frac{1}{t}\log U_{H,h}(t)&=&\max_{z\in F_t\setminus H_t}\biggl[\frac
{1}{t}\log U_{H,h}(t,z)\biggr]+o(1)\\
&\le&\max_{z\in F_t\setminus H_t}\biggl[
\Phi_t(z)+\frac{1}{t}\max_{h_t(z)\le n}\biggl\{\eta(n,z)-\eta(z)-
\frac{n-|z|}{2}\log\xi(z)\biggr\}\biggr]\\
&&{}+O(t^{q-\delta}).
\end{eqnarray*}
To prove (b), we show that the second term on the right-hand side of
\textup{(\ref{nadoelo})}
is negligible. Let $z\in F_t$. By Lemma \ref{asymp}(i) we have
$\xi(z)>t^{c}$ eventually.
Further, for $n\ge|z|$, we use Lemma \ref{nn} and the substitution
$r=n/|z|-1$ to get
%
\begin{eqnarray}\label{hd1}
&&\max_{h_t(z)\vee|z|\le n\le r_tg_t}\biggl\{\eta(n,z)-\eta(z)-\frac
{n-|z|}{2}\log\xi(z)\biggr\}\nonumber\\
&&\qquad\le\max_{n\ge|z|}\biggl[(n-|z|)\log\frac{2d
en}{(n-|z|)\sqrt{\xi(z)}}\biggr]+K\\
&&\qquad\le |z|\max_{r\ge0}\biggl[r\log\frac{2de
(r+1)}{rt^{c/2}}\biggr]+K.\nonumber
\end{eqnarray}
If $t$ is large enough, the expression in the square brackets is
negative for $r\ge1$,
hence the maximum is attained at some $r< 1$. Using this to estimate
the numerator and
optimizing the estimate, we obtain
%
\begin{eqnarray}\label{hd2}
\max_{r\ge0}\biggl[r\log\frac{2de(r+1)}{rt^{c/2}}\biggr]
\le\max_{r\ge0}\biggl[r\log\frac{4de}{rt^{c/2}}\biggr]
=4dt^{-c/2}.
\end{eqnarray}
Finally, since $|z|\le r_tg_t$, we obtain, combining \textup{(\ref{hd1})}
and \textup{(\ref{hd2})}
and, {if necessary,} decreasing $\delta$
so that it satisfies $\delta<c/2$,
\begin{eqnarray*}
\max_{h_t(z)\vee|z|\le n\le r_tg_t}\biggl\{\eta(n,z)-\eta(z)-\frac
{n-|z|}{2}\log\xi(z)\biggr\}
&\le&r_tg_t 4dt^{-c/2}+K\\
&=& O(t^{q+1-\delta}).
\end{eqnarray*}
Using this on the right-hand side of (a) completes the proof.
\end{pf}

\subsection{A lower bound for the growth of the mass}\label{sec-LB}

{We now derive a \textit{lower} bound for $U(t)$, which is a slight improvement
on the bound given in \cite{HMS08}, Lemma 2.2. This argument does not
rely on Lemma \ref{l10}.}

\begin{prop}
\label{iub}
Almost surely, eventually for all $t$
\begin{eqnarray*}
\frac{1}{t}\log U(t)\ge\Phi_t\bigl(Z_t^{({1})}\bigr)-2d+o(1).
\end{eqnarray*}
\end{prop}

\begin{pf} {The proof follows the}
same lines as in \cite{HMS08}, Lemma 2.2, so that we will shorten some
computations if they are the same.
Let ${\rho\in(0,1]}$ and $z\in\mathbb Z^d$ {with $|z|\ge2$}.
Denote by
\begin{eqnarray*}
A_{t}^{z,\rho}=\{J_{\rho t}=|z|, X_s=z\ \forall s\in[\rho
t,t]\}
\end{eqnarray*}
the event that the random walk $X$ reaches the point $z$ before time
$\rho t$, making the minimal
possible number of jumps, and stays at $z$ for the rest of the time.
Denote by $P_{\lambda}(\cdot)$ the Poisson distribution with
parameter $\lambda$.
Then, using Stirling's formula, we obtain
\begin{eqnarray*}
\mathbb P_0 (A_{t}^{z,\rho}) &=&
\frac{N(z)P_{2d\rho t}(|z|)P_{2d(1-\rho)t}(0)}{(2d)^{|z|}}\\
&=& \exp\biggl\{\eta(z)-|z|\log\frac{|z|}{e\rho
t}-2dt+{O(\log|z|)} \biggr\},
\end{eqnarray*}
where the last error term is bounded by the multiple of $\log|z|$
with an absolute constant.
As $\xi(z)\ge0$ almost surely for all $z$, we obtain, for all $\rho$
and $z$ {as above},
\begin{eqnarray*}
U(t)& =&\mathbb E_0 \biggl[\exp\biggl\{\int_0^t\xi(X_s)\,ds\biggr\}
\biggr] \ge e^{t(1-\rho)\xi(z)}\mathbb P_0 (A^{z,\rho}_t)
\\
&\ge&\exp\biggl\{t(1-\rho)\xi(z)+\eta(z)-|z|\log\frac
{|z|}{e\rho t}
-2dt + {O(\log|z|)}\biggr\}.
\end{eqnarray*}
Since {$\log|z| = o(t)$ for $|z|\le t^{\beta}$} for any fixed
positive $\beta$,
this implies
%
\begin{eqnarray}
\label{2max1}
\frac1t\log U(t)&\ge&\max_{0<\rho\le1}\max_{1\le|z|\le t^{\beta
}}\biggl[(1-\rho)\xi(z)+\frac{\eta(z)}{t}
-\frac{|z|}{t}\log\frac{|z|}{e\rho t}\biggr]\nonumber\\[-8pt]\\[-8pt]
&&{}-2d+o(1).\nonumber
\end{eqnarray}
Let $\widehat\eta\in(\frac{d}{\alpha},1)$ and $\beta=(1-\widehat
\eta)^{-1}(1+\varepsilon)$, $\varepsilon>0$.
By \cite{HMS08}, Lemma 3.5, there is $r_0$ such that $\xi_r^{({1})}\le
r^{\widehat\eta}$
for all $r>r_0$. We {thus} have, using the bound $\eta(z)\le|z|\log
d$ and {a}
similar computation as in \cite{HMS08}, Lemma 2.2,
%
\begin{eqnarray}\label{rem_z}
&&\max_{|z|>\max\{r_0,t^{\beta}\}}\biggl[(1-\rho)\xi(z)+\frac{\eta
(z)}{t}-\frac{|z|}{t}\log\frac{|z|}{e\rho t}\biggr]\nonumber\\
&&\qquad\le\max_{|z|>\max\{r_0,t^{\beta}\}}\biggl[(1-\rho)\xi
_{|z|}^{({1})}-\frac{|z|}{t}\log\frac{|z|}{de\rho t}\biggr]
\\
&&\qquad \le\max_{|z|>\max\{r_0,t^{\beta}\}}\biggl[|z|^{\widehat
\eta}\biggl(1-\rho-t^{\varepsilon}
\log\frac{t^{\beta-1}}{de\rho}\biggr)\biggr]<0,\nonumber
\end{eqnarray}
eventually for all $t$. Recall that {$\frac1t \log U(t)\ge0$}
and take $t$ large enough so that $t^{\beta}>r_0$.
Then \textup{(\ref{rem_z})} implies that the maximum in \textup
{(\ref{2max1})} can be
taken over all $z$
instead of $|z|\le t^{\beta}$. It is easy to see that this maximum is
attained at $\rho=\frac{|z|}{t\xi(z)}$
unless {this value exceeds one}. Substituting this $\rho$ into \textup
{(\ref{2max1})}
we obtain
\begin{eqnarray*}
\frac1t\log U(t)&\ge&\max_{z\in\mathbb Z^d}\biggl[\xi(z)+\frac
{\eta(z)}{t}
-\frac{|z|}{t}\log\xi(z)\biggr]{\mathbh1}\{t\xi(z)>|z|\}-2d+o(1)\\
&=&\Phi_t\bigl(Z_t^{({1})}\bigr)-2d+o(1),
\end{eqnarray*}
which completes the proof.
\end{pf}

\subsection{Negligible parts of the total mass}\label{sec-U12345}

In this section we show that the main contribution to the Feynman--Kac
formula for $U(t)$ comes from those
paths that pass through $Z_t^{({1})}$ or $Z_t^{({2})}$ and do not make
{significantly} more than $|Z_t^{({1})}|\wedge|Z_t^{({2})}|$ steps.
{For this purpose,} we define five path classes and show that the
latter four of them each give
a negligible contribution to the total mass $U(t)$.

In the sequel, we assume that $\delta$ is taken small enough {so that
Lemma \ref{l10} holds and $\delta<q$.
We decompose the set of all paths $[0,t]\to\mathbb Z^d$ into the following
five classes}:
\begin{eqnarray*}
A_i=\cases{
\biggl\{J_t\le r_tg_t, \exists
z\in\bigl\{Z_t^{({1})},Z_t^{({2})}\bigr\}\dvtx\cr
\hspace*{7.01pt}\displaystyle\max_{s\in[0,t]}\xi(X_s)=\xi(z), J_t< |z|(1+t^{-\delta
/2})\biggr\},&\quad $i=1$,\cr
\biggl\{J_t\le r_tg_t, \exists
z\in\bigl\{Z_t^{({1})},Z_t^{({2})}\bigr\}\dvtx\cr
\hspace*{7.01pt}\displaystyle\max_{s\in[0,t]}\xi(X_s)=\xi(z), J_t\ge|z|(1+t^{-\delta
/2})\biggr\},&\quad $i=2$,\cr
\biggl\{J_t\le r_tg_t, \exists z\in F_t\setminus\bigl\{
Z_t^{({1})},Z_t^{({2})}\bigr\}\dvtx
\displaystyle\max_{s\in[0,t]}\xi(X_s)=\xi(z)\biggr\},&\quad $i=3$,\cr
\biggl\{J_t\le r_tg_t, \displaystyle\max_{s\in[0,t]}\xi(X_s)\le\xi
_{r_tg_t}^{({k_t})}\biggr\},&\quad $i=4$,\cr
\{J_t>r_tg_t\},&\quad $i=5$
}
\end{eqnarray*}
and split the total mass into five components $U(t)=\sum_{i=1}^5
U_i(t)$, where
\begin{eqnarray*}
U_i(t)=\mathbb E_0\biggl[\exp\biggl\{\int_0^t\xi(X_s)\,ds\biggr\}{\mathbh1}
_{A_i}\biggr],\qquad 1\le i\le5.
\end{eqnarray*}

\begin{lemma}
\label{25}
Almost surely, $\lim_{t\to\infty}U_i(t)/U(t)=0$ for $2\le i\le5$.
\end{lemma}

\begin{pf} \textit{Case} $i=2$: Denote $h_t(z)=|z|(1+t^{-\delta/2})$
and $H_t=F_t\setminus\{
Z_t^{({1})},Z_t^{({2})}
\}$.
By Lemma \ref{l10}(a),
%
\begin{eqnarray}\label{n25}
&&\frac{1}{t}\log U_2(t)\nonumber\\
&&\qquad\le \frac{1}{t}\log
U_{H,h}(t)\nonumber\\[-8pt]\\[-8pt]
&&\qquad\le\max_{z\in\{Z_t^{({1})},Z_t^{({2})}\}}\biggl\{\Phi_t(z)+\frac{1}{t}
\max_{{n}\ge|z|(1+t^{-\delta})}\biggl[\eta(n,z)-\eta(z)\nonumber\\
&&\qquad\quad\hspace*{153.4pt}{}-\frac
{n-|z|}{2}\log\xi(z)\biggr]\biggr\} + O(t^{q-\delta}).\hspace*{-25pt}\nonumber
\end{eqnarray}
For each $z\in\{Z_t^{({1})},Z_t^{({2})}\}$ we have by Lemma \ref
{mon}(i), for any
$c>0$, that $\xi(z)>(2de)^2t^{q-c}$ eventually.
Together with Lemma \ref{nn} this implies
\begin{eqnarray*}
&&\max_{n\ge|z|(1+t^{-\delta/2})}\biggl[\eta(n,z)-\eta(z)-\frac
{n-|z|}{2}\log\xi(z)\biggr]\\
&&\qquad\le\max_{n\ge|z|(1+t^{-\delta/2})}\biggl[(n-|z|)\log\frac
{2den\xi(z)^{-{1}/{2}}}{n-|z|}\biggr]+K\\
&&\qquad\le\max_{n\ge|z|(1+t^{-\delta/2})}\biggl[(n-|z|)\log\frac
{nt^{-({q-c})/{2}}}{n-|z|}\biggr]+K\\
&&\qquad=|z|\max_{r\ge t^{-\delta/2}}\biggl[r\log\frac{(r+1)t^{-
({q-c})/{2}}}{r}\biggr]+K.
\end{eqnarray*}
It is easy to check that, eventually for all $t$, the function under
the maximum is decreasing on $(t^{-({q-c})/{2}},\infty)$ if $c<q$.
Since $\delta<q$ we can choose $c$ so small that
$\delta<q-c$. The maximum is then attained at $r=t^{-\delta/2}$ and,
as 
$|z|\ge t^{q+1-{\delta}/{4}}$ by Lemma \ref{mon}(iii), we obtain
\begin{eqnarray*}
\max_{n\ge|z|(1+t^{-\delta/2})}\biggl[\eta(n,z)-\eta(z)+\frac
{n-|z|}{2}\log\xi(z)\biggr]
&\le&-|z| t^{-\delta/2}\log\bigl(t^{({q-c-\delta})/{2}}\bigr)+K\\
&\le&-t^{q+1-{3\delta}/{4}}
\end{eqnarray*}
eventually for all $t$. Combining this with \textup{(\ref{n25})} and using
Proposition \ref{iub} we finally get
\begin{eqnarray*}
\frac{1}{t}\log\frac{U_{2}(t)}{U(t)}
&\le&\max_{i=1,2}\bigl\{\Phi_t\bigl(Z_t^{({i})}\bigr)-t^{q-{3\delta
}/{4}}\bigr\}-\Phi_t\bigl(Z_t^{({1})}\bigr)+O(t^{q-\delta})\\
&=&-t^{q-{3\delta}/{4}}+O(t^{q-\delta})\to-\infty.
\end{eqnarray*}

\textit{Case} $i=3$: Pick $h_t(z)=0$ and $H_t=\{
Z_t^{({1})},Z_t^{({2})}\}$.
By Lemma \ref{l10}(b) we obtain
\begin{eqnarray*}
\frac{1}{t}\log U_3(t)&\le&\frac{1}{t}\log U_{H,h}(t)\le
\max_{z\in\mathbb Z^d\setminus\{Z_t^{({1})},Z_t^{({2})}\}}\Phi
_t(z)+O(t^{q-\delta})\\
&=&\Phi_t\bigl(Z_t^{({3})}\bigr)+O(t^{q-\delta}).
\end{eqnarray*}
It remains to apply Propositions \ref{iub} and \ref{gap} to get eventually
\begin{eqnarray*}
\frac{1}{t}\log\frac{U_3(t)}{U(t)}\le\Phi_t\bigl(Z_t^{({3})}\bigr)-\Phi
_t\bigl(Z_t^{({1})}\bigr)+O(t^{q-\delta})
\le-a_t\lambda_t+O(t^{q-\delta})\to-\infty.
\end{eqnarray*}

\textit{Case} $i=4$: We estimate the integral in the Feynman--Kac
formula by $t\xi
_{r_tg_t}^{({k_t})}$.
Lemma \ref{asymp}(i) implies that there is $c>0$ such that eventually
\begin{eqnarray*}
\frac{1}{t}\log U_4(t)\le\xi_{r_tg_t}^{({k_t})}\le t^{q-c}.
\end{eqnarray*}
On the other hand, it follows from \cite{HMS08}, Theorem 1.1, that, for
each $\widehat c>0$,
we have $\frac{1}{t}\log U(t)\ge t^{q-\widehat c}$
eventually. Since $\widehat c$ can be taken smaller than $c$, the
statement is proved.

\textit{Case} $i=5$: We decompose the Feynman--Kac formula according
to the number~$J_t$ of jumps.
Observe that the integral there can be estimated by $t\xi
_{J_t}^{({1})}$ and use that
$J_t$ has Poisson distribution with parameter $2dt$. Thus, we obtain
\begin{eqnarray*}
U_5(t)&=& \sum_{n>r_tg_t} \mathbb E_0\biggl[\exp\biggl\{\int_0^t\xi(X_s)\,d
s\biggr\}{\mathbh1}\{J_t=n\}\biggr]\\
&\le&\sum_{n>r_tg_t} \exp\biggl\{t\xi_n^{({1})}-2dt+\log\frac
{(2dt)^n}{n!}\biggr\}.
\end{eqnarray*}
Pick $0<\varepsilon<\nu/(q+1)$ and assume that $t$ is large enough.
It follows from \cite{HMS08}, Lemma 3.5, that
$\xi_n^{({1})}<n^{{d}/{\alpha}}(\log n)^{{1}/{\alpha
}+\varepsilon}$ for all $n>r_tg_t$.
Further, it follows from Stirling's formula that $n!>(n/e)^n$ for
all $n>r_tg_t$.
Then, for all $n>r_tg_t$, we obtain, using monotonicity in $n$,
\[
t\xi_n^{({1})}-2dt+\log\frac{(2dt)^n}{n!}
< tn^{{d}/{\alpha}}(\log n)^{{1}/{\alpha}+\varepsilon
}+n\log\frac{2det}{n}\le-n^{{d}/{\alpha}}.
\]
Combining the last two displays we obtain that $U_5(t)=o(1)$.
\end{pf}

\subsection{An upper bound for the growth of the mass}\label{sec-UB}

Lemmas \ref{l10} and \ref{25} make it possible to prove an upper bound
for $\frac1t\log U(t)$, which is asymptotically
equal to the lower bound of Proposition \ref{iub}.

\begin{prop}\label{cor-uppboundU}
Fix $\delta>0$ as in Lemma \ref{l10}. Then, almost surely, eventually
for all $t$,
\begin{eqnarray*}
\frac{1}{t}\log U(t)\le\Phi_t\bigl(Z_t^{({1})}\bigr)+O(t^{q-\delta}).
\end{eqnarray*}
\end{prop}

\begin{pf} Consider $H_t=\varnothing$ and $h_t=0$. Then
$U_{H,h}(t)=U_1(t)+U_2(t)+U_3(t)$.
Since for the remaining two functions we have
$U_4(t)+U_5(t)\le U(t)o({1})$ by Lemma \ref{25}, we obtain by
Lemma \ref{l10}(b) that
$\frac{1}{t}\log U(t)\le\frac{1}{t}\log U_{H,h}(t)(1+o(1))\le\Phi
_t(Z_t^{({1})})
+O(t^{q-\delta})$.
\end{pf}

\section{Almost sure localization in two points}\label{s_as}

In this section, we prove Theorem~\ref{main_as}. In Section \ref
{sec-threeevents} we introduce a decomposition
into three events, formulate our main steps and provide some technical
preparation. The remaining
Sections~\mbox{\ref{sec-Proof1}--\ref{sec-Proof3}} give the proofs of the
localization on the three respective events.

\subsection{Decomposition into three events}\label{sec-threeevents}

In the proof of Theorem \ref{main_as}, we distinguish between three
disjoint events
{constituting a partition of} the full probability space:
\begin{itemize}
\item$\Phi_t(Z_t^{({1})})-\Phi_t(Z_t^{({2})})$ is small and the
sites $Z_t^{({1})}$ and $Z_t^{({2})}$ are close to
each other.
\item$\Phi_t(Z_t^{({1})})-\Phi_t(Z_t^{({2})})$ is small but the
sites $Z_t^{({1})}$ and $Z_t^{({2})}$ are far away
from each other.
\item$\Phi_t(Z_t^{({1})})-\Phi_t(Z_t^{({2})})$ is large.
\end{itemize}
We prove the two point localization on each event {by different arguments.}
To be precise, for $i=1,2$, denote by
\begin{eqnarray*}
\Gamma_t^{({i})}=\bigl\{z\in\mathbb Z^d\dvtx\bigl|z-Z_t^{({i})}\bigr|+\min\bigl\{
|z|,\bigl|Z_t^{({i})}\bigr|\bigr\}
< \bigl|Z_t^{({i})}\bigr|(1+t^{-\delta/2})\bigr\}
\end{eqnarray*}
the set containing all sites $z$ such that there is a path
of length less than $|Z_t^{({i})}|(1+t^{-\delta/2})$ starting from
the origin passing through both $z$ and $Z_t^{({i})}$.
Further, denote
\begin{eqnarray*}
\Gamma_t=\bigl\{z\in\mathbb Z^d\dvtx\bigl|z-Z_t^{({1})}\bigr|+\min\bigl\{
|z|,\bigl|Z_t^{({1})}\bigr|\bigr\}< \bigl|Z_t^{({1})}
\bigr|(1+6t^{-\delta/2})\bigr\}.
\end{eqnarray*}
In Sections \ref{sec-Proof1}--\ref{sec-Proof3} we prove the following
three propositions:

\begin{prop}
\label{ttwo}
Almost surely,
\begin{eqnarray*}
&&\lim_{t\to\infty}
\Biggl[U(t)^{-1}
\sum_{z\in\mathbb Z^d\setminus\{Z_t^{({1})},Z_t^{({2})}\}}
u(t,z)\Biggr]\\
&&\qquad{}\times {\mathbh1}\bigl\{\Phi_t\bigl(Z_t^{({1})}\bigr)-\Phi_t\bigl(Z_t^{({2})}\bigr)< a_t\lambda
_t/2, Z_t^{({2})}\in\Gamma_t^{({1})}\bigr\}=0.
\end{eqnarray*}
\end{prop}

\begin{prop}
\label{two}
Almost surely,
\begin{eqnarray*}
&&\lim_{t\to\infty}
\Biggl[U(t)^{-1}
\sum_{z\in\mathbb Z^d\setminus\{Z_t^{({1})},Z_t^{({2})}\}}
u(t,z)\Biggr]\\
&&\qquad{}\times {\mathbh1}\bigl\{\Phi_t\bigl(Z_t^{({1})}\bigr)-\Phi_t\bigl(Z_t^{({2})}\bigr)< a_t\lambda
_t/2, Z_t^{({2})}\notin\Gamma_t^{({1})}\bigr\}=0.
\end{eqnarray*}
\end{prop}

\begin{prop}
\label{three}
Almost surely,
\begin{eqnarray*}
\lim_{t\to\infty}
\Biggl[U(t)^{-1} \sum_{z\in\mathbb Z^d\setminus
\{Z_t^{({1})}
\}} u(t,z)\Biggr]
{\mathbh1}\bigl\{\Phi_t\bigl(Z_t^{({1})}\bigr)-\Phi_t\bigl(Z_t^{({2})}\bigr)\ge a_t\lambda
_t/2\bigr\}=0.
\end{eqnarray*}
\end{prop}

Obviously, Theorem \ref{main_as} follows immediately from the three
propositions.
For each of them, the idea of the proof is to decompose $u$ into a sum
of two functions
$u_1$ and $u_2$ (which is different {in different cases}) such that
$u_2$ is negligible
and localization of $u_1$ can be shown with the help of our spectral
bounds derived in
Section \ref{sec-spectral}. If the gap between $\Phi_t(Z_t^{({1})})$
and $\Phi_t(Z_t^{({2})})$ is
small (Cases 1 and 2)
then both points $Z_t^{({1})}$ and $Z_t^{({2})}$ contribute to the
total mass.
However, the strategy of
the proof is different, since in the second case the points
$Z_t^{({1})}$ and $Z_t^{({2})}$ do not interact as they are far away
from each other,
whereas in the first case
they do. If the gap between $\Phi_t(Z_t^{({1})})$ and $\Phi
_t(Z_t^{({2})})$ is large (Case 3) only
the site $Z_t^{({1})}$ contributes to
the total mass. In the remaining part of this section, we prove a
lemma, which {is
used in the proof of} each of the three propositions.

\begin{lemma}
\label{mb}
There is ${c\in(0,q)}$ such that, almost surely eventually for all $t$:
\begin{longlist}
\item $\xi(z)< \xi(Z_t^{({1})})-t^{q-c}$ for all $z\in\Gamma
_t\setminus\{Z_t^{({1})},Z_t^{({2})}\}$,
\item  $\xi(z)< \xi(Z_t^{({1})})-t^{q-c}$ for all $z\in
\Gamma
_t\setminus\{Z_t^{({1})}\}$
if $\Phi_t(Z_t^{({1})})-\Phi_t(Z_t^{({2})})\ge a_t\lambda_t/2$,
\item  $\xi(z)< \xi(Z_t^{({2})})-t^{q-c}$ for all $z\in
\Gamma
_t^{({2})}\setminus\{Z_t^{({1})},Z_t^{({2})}\}$
if $\Phi_t(Z_t^{({1})})-\Phi_t(Z_t^{({2})})< a_t\lambda_t/2$,
\item  $\Gamma_t^{({1})}\subset\Gamma_t$. If $Z_t^{({2})}\in
\Gamma_t^{({1})}$
then $\Gamma_t^{({2})}\subset\Gamma_t$.
\end{longlist}
\end{lemma}

\begin{pf}
We prove (i)--(iii) simultaneously, first making the following
observations:
\begin{longlist}[(3)]
\item[(1)] By Proposition \ref{gap} we have $\Phi
_t(Z_t^{({1})})-\Phi
_t(z)>a_t\lambda_t/2$ for all $z\notin\{Z_t^{({1})},Z_t^{({2})}\}$.
\item[(2)] $\Phi_t(Z_t^{({1})})-\Phi_t(z)\ge\Phi
_t(Z_t^{({1})})-\Phi_t(Z_t^{({2})})\ge a_t\lambda_t/2$ for
all $z\neq Z_t^{({1})}$ by assumption.
\item[(3)] Using Proposition \ref{gap} and our assumption we
obtain $\Phi_t(Z_t^{({2})})-\Phi_t(Z_t^{({3})})>a_t\lambda_t/2$.
Hence $\Phi_t(Z_t^{({2})})-\Phi_t(z)\ge\Phi_t(Z_t^{({2})})-\Phi
_t(Z_t^{({3})})>a_t\lambda_t/2$ for all $z\notin\{
Z_t^{({1})},Z_t^{({2})}\}$.
\end{longlist}
Thus, to show (i)--(iii), it suffices to prove that there exists
${c\in(0,q)}$
such that eventually
\[
\xi\bigl(Z_t^{({i})}\bigr)-\xi(z)>t^{q-c}
\]
for each $i\in\{1,2\}$
and each $z$ satisfying
%
\begin{eqnarray}
\label{iiii}
\Phi_t\bigl(Z_t^{({i})}\bigr)-\Phi_t(z)&\ge& a_t\lambda_t/2
\quad\mbox{and}\nonumber\\[-8pt]\\[-8pt]
\bigl|z-Z_t^{({i})}\bigr|+\min\bigl\{|z|,\bigl|Z_t^{({i})}\bigr|\bigr\}
&<&\bigl|Z_t^{({i})}\bigr|(1+6t^{-\delta/2}).\nonumber
\end{eqnarray}
Assume that the statement is false. {Then given $c<q$ there is an
arbitrarily large $t$ and $z\in\mathbb Z^d$ satisfying \textup{(\ref
{iiii})} with
$\xi(Z_t^{({i})})-\xi(z)
\le t^{q-c}$. Then}
%
\begin{eqnarray}
\label{cc1}
\Phi_t\bigl(Z_t^{({i})}\bigr)-\Phi_t(z)
&=&\bigl[\xi\bigl(Z_t^{({i})}\bigr)-\xi(z)\bigr]+\frac{|Z_t^{({i})}|}{t}\log
\frac{\xi(z)}{\xi(Z_t^{({i})})}\nonumber\\[-8pt]\\[-8pt]
&&{} +\frac{|z|-|Z_t^{({i})}|}{t}\log\xi(z)+\frac
{\eta(Z_t^{({1})})-\eta(z)}{t}.\nonumber
\end{eqnarray}
We can bound the second summand by zero if $\xi(z)\le\xi
(Z_t^{({i})})$. For $\xi(z)> \xi(Z_t^{({i})})$,
we use the inequality $\log x\le x-1$ for $x>0$ to obtain, by
Lemma \ref{mon}(ii), eventually
\begin{eqnarray*}
\frac{|Z_t^{({i})}|}{t}\log\frac{\xi(z)}{\xi(Z_t^{({i})})}\le
\frac{|Z_t^{({i})}|(\xi(z)-\xi(Z_t^{({i})}))}{t\xi
(Z_t^{({i})})}<\xi(z)-\xi\bigl(Z_t^{({i})}\bigr).
\end{eqnarray*}
In both cases we obtain the estimate for the first two terms
\begin{eqnarray*}
\bigl[\xi\bigl(Z_t^{({i})}\bigr)-\xi(z)\bigr]+\frac{|Z_t^{({i})}|}{t}\log
\frac{\xi(z)}{\xi(Z_t^{({i})})}
\le\max\bigl\{\xi\bigl(Z_t^{({i})}\bigr)-\xi(z),0\bigr\}\le t^{q-c}=o(a_t\lambda_t).
\end{eqnarray*}
We prove that the remaining two terms {in \textup{(\ref{cc1})}} are
of order
$o(a_t\lambda_t)$
as well. First, assume that $|z|<|Z_t^{({i})}|$. Then \textup{(\ref
{iiii})} implies
%
\begin{eqnarray}
\label{nl1}
\bigl|Z_t^{({i})}\bigr|\le\bigl|z-Z_t^{({i})}\bigr|+|z|<\bigl|Z_t^{({i})}\bigr|(1+6t^{-\delta/2}).
\end{eqnarray}
Notice that $\eta(Z_t^{({i})})\le\eta
(|z-Z_t^{({i})}|+|Z_t^{({i})}|,z)$ as
to any path of length $|Z_t^{({i})}|$ passing through $Z_t^{({i})}$
we can add a path of length $|z-Z_t^{({i})}|$ in such a way that it passes
through $z$. Using Lemmas \ref{nn} and \ref{mon}(iii)
and \textup{(\ref{nl1})} we obtain
\begin{eqnarray*}
&&\eta\bigl(Z_t^{({1})}\bigr)-\eta(z)\\
&&\qquad\le\eta\bigl(\bigl|z-Z_t^{({i})}\bigr|+\bigl|Z_t^{({i})}\bigr|,z\bigr)-\eta(z)\\
&&\qquad\le\bigl(\bigl|z-Z_t^{({i})}\bigr|+\bigl|Z_t^{({i})}\bigr|-|z|\bigr)\log\frac{2d
{e}(|z-Z_t^{({i})}|+|Z_t^{({i})}|)}{|z-Z_t^{({i})}|+|Z_t^{({i})}|-|z|}+K\\
&&\qquad\le\bigl(6t^{-\delta/2}\bigl|Z_t^{({i})}\bigr|+2\bigl(\bigl|Z_t^{({i})}\bigr|-|z|\bigr)\bigr)
\log\frac{2de((2+t^{-\delta
/2})|Z_t^{({i})}|-|z|)}{|Z_t^{({i})}|-|z|}+K\\
&&\qquad\le2\bigl(\bigl|Z_t^{({i})}\bigr|-|z|\bigr)
\log\frac{5d
{e}|Z_t^{({i})}|}{|Z_t^{({i})}|-|z|}+O(t^{q+1+\varepsilon-\delta/4}).
\end{eqnarray*}
Substituting this as well as the estimate for the first two terms
into \textup{(\ref{cc1})} we obtain
\begin{eqnarray*}
\Phi_t\bigl(Z_t^{({i})}\bigr)-\Phi_t(z)\le\frac{2(|Z_t^{({i})}|-|z|)}{t}
\log\frac{5de|Z_t^{({i})}|}{(|Z_t^{({i})}|-|z|)\sqrt{\xi
(z)}}+o(a_t\lambda_t).
\end{eqnarray*}
By Lemma \ref{mon}(i) we have $\xi(Z_t^{({i})})>t^{q-c/4}$ as
$t$ is large enough. By
assumption we then have $\xi(z)\ge\xi
(Z_t^{({i})})-t^{q-c}>t^{q-c/2}$. Hence the expression
under the logarithm is positive only if $|Z_t^{({i})}|-|z|<5d
e|Z_t^{({i})}| t^{-q/2+c/4}$,
which is smaller than $t^{q/2+1+c/2}$ by Lemma \ref{mon}(iii).
Since $c<q$ we obtain $\Phi_t(Z_t^{({i})})-\Phi_t(z)\le
o(a_t\lambda_t)$.

Finally, consider $|z|\ge|Z_t^{({i})}|$.
For the third term in \textup{(\ref{cc1})} we notice that \textup
{(\ref{iiii})} implies that
$|z|\le|Z_t^{({i})}|(1+6t^{-\delta/2})$. Then we use Lemma \ref
{mon}(iii) and
\cite{HMS08}, Lemma~3.5, which gives
\begin{eqnarray*}
\frac{|z|-|Z_t^{({i})}|}{t}\log\xi(z)\le6t^{-\delta
/2-1}\bigl|Z_t^{({i})}\bigr|
\log\bigl(\bigl|Z_t^{({i})}\bigr|^{{d}/{\alpha}+\delta}(1+6t^{-\delta
/2})^{{d}/{\alpha}+\delta}\bigr)
=o(a_t\lambda_t).
\end{eqnarray*}
For the last term in \textup{(\ref{cc1})} we obtain from \textup
{(\ref{iiii})} that
$\eta(Z_t^{({i})})\le\eta(|Z_t^{({i})}|(1+6t^{-\delta/2}),z)$.
Hence, by Lemma \ref{nn},
\begin{eqnarray*}
\eta\bigl(Z_t^{({1})}\bigr)-\eta(z)
&\le&\eta\bigl(\bigl|Z_t^{({i})}\bigr|(1+6t^{-\delta/2}),z\bigr)-\eta(z)\\
&\le&\bigl(\bigl|Z_t^{({i})}\bigr|(1+6t^{-\delta/2})-|z|\bigr)\log
\frac{2de|Z_t^{({i})}|(1+6t^{-\delta
/2})}{|Z_t^{({i})}|(1+6t^{-\delta/2})-|z|}+K.
\end{eqnarray*}
Notice that, for $a>0$, $x\mapsto x\log\frac{a}{x}$ is an increasing
function on $(0,a/e)$.
Since $|z|\ge|Z_t^{({i})}|$ we have $|Z_t^{({i})}|(1+6t^{-\delta
/2})-|z|\le6t^{-\delta/2}|Z_t^{({i})}|$,
which is smaller than $2d|Z_t^{({i})}|(1+6t^{-\delta/2})$. Hence we obtain
\begin{eqnarray*}
\eta\bigl(Z_t^{({1})}\bigr)-\eta(z)\le6t^{-\delta/2}\bigl|Z_t^{({i})}\bigr|\log\frac
{de(1+6t^{-\delta/2})}{3t^{-\delta/2}}\le O(t^{q+1-\delta/4}).
\end{eqnarray*}
This proves that the last term in \textup{(\ref{cc1})} is also
bounded by
$o(a_t\lambda_t)$.
It remains to notice that we have proved $\Phi_t(Z_t^{({i})})-\Phi
_t(z)\le o(a_t\lambda_t)$,
which contradicts our assumption that $\Phi_t(Z_t^{({i})})-\Phi
_t(z)\ge a_t\lambda_t{/2}$.
This proves (i)--(iii).

(iv) The first statement is trivial. To prove the second one, we
pick ${z\in\Gamma_t^{({2})}}$. For any such $z$ there exists a path
of length
less than $|Z_t^{({2})}|(1+t^{\delta/2})$ starting at the origin and going
through $z$ and $Z_t^{({2})}$. If $|z|\le|Z_t^{({2})}|$
we can choose the path in such a way that it ends at $Z_t^{({2})}$. If
$|z|>|Z_t^{({2})}|$ then
$|z-Z_t^{({2})}|<|Z_t^{({2})}| t^{\delta/2}$ and so there is a path
of length less
than $|Z_t^{({2})}|(1+2t^{\delta/2})$
starting at the origin, going through~$z$ and ending in $Z_t^{({2})}$. In
either of the cases,
there is then a path of length less than $|Z_t^{({2})}|(1+2t^{\delta
/2})+|Z_t^{({1})}
-Z_t^{({2})}|$ starting
in the origin, passing through $z$, $Z_t^{({2})}$ and ending at $Z_t^{({1})}$.

Observe that since
$Z_t^{({2})}\in\Gamma_t^{({1})}$, there is a path of length less than
$|Z_t^{({1})}|(1+t^{-\delta/2})$
going through $Z_t^{({1})}$ and $Z_t^{({2})}$, which in particular
implies $|Z_t^{({2})}
|<|Z_t^{({1})}|(1+t^{-\delta/2})$.
If $|Z_t^{({2})}|<|Z_t^{({1})}|$, then $Z_t^{({2})}\in\Gamma
_t^{({1})}$ implies $|Z_t^{({2})}
|+|Z_t^{({1})}-Z_t^{({2})}|<|Z_t^{({1})}|(1+t^{-\delta/2})$ and so
\begin{eqnarray*}
\bigl|Z_t^{({2})}\bigr|(1+2t^{\delta/2})+\bigl|Z_t^{({1})}-Z_t^{({2})}\bigr|\le
\bigl|Z_t^{({1})}\bigr|(1+3t^{-\delta/2})<\bigl|Z_t^{({1})}
\bigr|(1+5t^{-\delta/2}).
\end{eqnarray*}
If $|Z_t^{({2})}|\ge|Z_t^{({1})}|$, then $Z_t^{({2})}\in\Gamma
_t^{({1})}$ implies
$|Z_t^{({1})}-Z_t^{({2})}|<|Z_t^{({1})}|t^{-\delta/2}$ and so
\begin{eqnarray*}
\bigl|Z_t^{({2})}\bigr|(1+2t^{\delta/2})+\bigl|Z_t^{({1})}-Z_t^{({2})}\bigr|
&\le&\bigl|Z_t^{({1})}\bigr|(1+t^{-\delta/2})(1+2t^{-\delta
/2})+\bigl|Z_t^{({1})}\bigr|t^{-\delta/2}\\
&<&\bigl|Z_t^{({1})}\bigr|(1+5t^{-\delta/2}).
\end{eqnarray*}
In each case {we obtain} that $z\in\Gamma_t$, {which completes the proof.}
\end{pf}

\subsection{First event: $\Phi_t(Z_t^{({1})})$ is close
to $\Phi_t(Z_t^{({2})})$ and $Z_t^{({1})}$ is close to
$Z_t^{({2})}$}\label{sec-Proof1}

In this section we prove Proposition \ref{ttwo}.
Let us decompose $u$ into $u=u_1+u_2$ with
\begin{eqnarray*}
u_1(t,z)&=&\mathbb E_0\biggl[\exp\biggl\{\int_0^t\xi(X_s)\,ds\biggr\}
{\mathbh1}\{
X_t=z\}
{\mathbh1}\bigl\{\tau_{\{Z_t^{({1})},Z_t^{({2})}\}}\le t,\tau_{\Gamma
_t^{\mathrm{c}}}>t\bigr\}
\biggr],\\
u_2(t,z)&=&\mathbb E_0\biggl[\exp\biggl\{\int_0^t\xi(X_s)\,ds\biggr\}
{\mathbh1}\{
X_t=z\}
{\mathbh1}\bigl\{\tau_{\{Z_t^{({1})},Z_t^{({2})}\}}> t\mbox{ or }\tau
_{\Gamma_t^{\mathrm
{c}}}\le t\bigr\}\biggr].
\end{eqnarray*}
We show that $u_1$ is localized in $Z_t^{({1})}$ and $Z_t^{({2})}$ (see
Lemma \ref
{lem-u1esti1})
and that the contribution of $u_2$ is negligible (see Lemma \ref{lem-u1esti2}).

\begin{lemma}\label{lem-u1esti1} Almost surely,
\begin{eqnarray*}
&&\lim_{t\to\infty}
\Biggl[U(t)^{-1}
\sum_{z\in\mathbb Z^d\setminus\{Z_t^{({1})},Z_t^{({2})}\}}
u_1(t,z)\Biggr]\\
&&\qquad{}\times {\mathbh1}\bigl\{\Phi_t\bigl(Z_t^{({1})}\bigr)-\Phi_t\bigl(Z_t^{({2})}\bigr)< a_t\lambda
_t/2, Z_t^{({2})}\in\Gamma_t^{({1})}\bigr\}=0.
\end{eqnarray*}
\end{lemma}

\begin{pf} We further split $u_1$ into three contributions $u_1=
u_{1,1}+u_{1,2}+u_{1,3}$ with
\begin{eqnarray*}
u_{1,j}(t,z)=\mathbb E_0\biggl[\exp\biggl\{\int_0^t\xi(X_s)\,ds\biggr\}
{\mathbh1}
\{X_t=z\}{\mathbh1}_{C_{1j}}\biggr]
\end{eqnarray*}
with
\begin{eqnarray*}
C_{1j}=\cases{
\bigl\{\tau_{Z_t^{({1})}}\le t,\tau_{\Gamma_t^{\mathrm{c}}}>t\bigr\}
, &\quad
$j=1$,\cr
\bigl\{\tau_{Z_t^{({1})}}> t,\tau_{Z_t^{({2})}}\le t,\tau_{[\Gamma
_t^{({2})}]^{\mathrm{c}}}>t\bigr\}, &\quad $j=2$,\cr
\bigl\{\tau_{Z_t^{({1})}}> t,\tau_{Z_t^{({2})}}\le t,\tau_{[\Gamma
_t^{({2})}]^{\mathrm{c}}}\le t,
\tau_{\Gamma_t^{\mathrm{c}}}>t\bigr\}, &\quad $j=3$.}
\end{eqnarray*}
Observe that the sets $C_{11}, C_{12},C_{13}$ are disjoint on the event
$\{Z_t^{({2})}\in\Gamma_t^{({1})}\}$ since $\Gamma_t^{({2})}\subset
\Gamma_t$
by Lemma \ref{mb}(iv). Furthermore, on this event, their union is
equal to the
event
\[
\bigl\{\tau_{\{Z_t^{({1})},Z_t^{({2})}\}}\le t, \tau_{\Gamma
_t^{\mathrm
{c}}}>t\bigr\}
\]
appearing in
the definition of $u_1(t,z)$. Hence, we indeed have
$u_1=u_{1,1}+u_{1,2}+u_{1,3}$.

We {now fix a sufficiently large $t$ and
argue on the event $\{\Phi_t(Z_t^{({1})})-\Phi_t(Z_t^{({2})})<
a_t\lambda_t/2, Z_t^{({2})}\in\Gamma
_t^{({1})}\}$.
We also fix some $c\in(0,q)$ and use this to distinguish between two cases.}

(1) First, we {assume} $\xi(Z_t^{({2})})\le\xi(Z_t^{({1})})-t^{q-c}$.
We show that $u_{1,1}$ and $u_{1,2}$ are localized around $Z_t^{({1})}$ and
$Z_t^{({2})}$,
respectively, and that the contribution of $u_{1,3}$ is negligible.

Let us fix $t$ large enough and pick
$B=\Gamma_t$, $\Omega=\{Z_t^{({1})}\}$ to study $u_{1,1}$
and $B=\Gamma_t^{({2})}\setminus\{Z_t^{({1})}\}$, $\Omega=\{
Z_t^{({2})}\}$ to
study $u_{1,2}$.
For the first {choice} we have
\begin{eqnarray*}
\mathfrak{g}_{\Omega,B}&=&\xi\bigl(Z_t^{({1})}\bigr)- \max_{\Gamma
_t\setminus\{
Z_t^{({1})}\}}\xi(z)\\
&=&\min\biggl\{\xi\bigl(Z_t^{({1})}\bigr)- \max_{\Gamma_t\setminus\{
Z_t^{({1})},Z_t^{({2})}
\}}\xi(z),
\xi\bigl(Z_t^{({1})}\bigr)-\xi\bigl(Z_t^{({2})}\bigr)\biggr\}\ge t^{q-c}
\end{eqnarray*}
by our assumption and Lemma \ref{mb}(i). For the second choice
we get
\begin{eqnarray*}
\mathfrak{g}_{\Omega,B}=\xi\bigl(Z_t^{({2})}\bigr)-\max_{z\in\Gamma
_t^{(2)}\setminus\{Z_t^{({1})},Z_t^{({2})}\}}\xi(z)\ge t^{q-c}
\end{eqnarray*}
by Lemma \ref{mb}(iii), using that $\Phi_t(Z_t^{({1})})-\Phi
_t(Z_t^{({2})})< a_t\lambda_t/2$.
Now we apply Lemma \ref{spectral} and use the monotonicity of $\varphi$
to obtain
%
\begin{eqnarray}
\label{2in}
\frac{\sum_{z\in\Gamma_t\setminus\{Z_t^{({1})}\}
}u_{1,1}(t,z)}{\sum
_{z\in\Gamma_t}u_{1,1}(t,z)}&\le&\varphi(t^{q-c})
\quad\mbox{and}\nonumber\\[-8pt]\\[-8pt]
\frac{\sum_{z\in\Gamma_t^{(2)}\setminus\{Z_t^{({1})},Z_t^{({2})}\}
}u_{1,2}(t,z)}{\sum_{z\in\Gamma_t^{({2})}\setminus\{Z_t^{({1})}\}}
u_{1,2}(t,z)}&\le&\varphi(t^{q-c}).\nonumber
\end{eqnarray}
Obviously, the estimate remains true if we increase the denominators
and sum over all
$z$ the larger function $u(t,z)$, which will produce $U(t)$. For the
numerators, notice that
$u_{1,1}(t,z)=0$ for all $z\notin\Gamma_t$ as the paths from
${C_{11}}$ do not leave
$\Gamma_t$, and $u_{1,2}(t,z)=0$ for all $z\notin\Gamma
_t^{({2})}\setminus\{Z_t^{({1})}\}$
as the paths from ${C_{12}}$ do not leave this set. Hence \textup
{(\ref{2in})} implies
\begin{eqnarray*}
U(t)^{-1} \sum_{z\in\mathbb Z^d\setminus\{
Z_t^{({1})}\}}
u_{1,1}(t,z)\le\varphi(t^{q-c})=o(1)
\end{eqnarray*}
and
\begin{eqnarray*}
U(t)^{-1} \sum_{z\in\mathbb Z^d\setminus\{Z_t^{({2})}\}
}
u_{1,2}(t,z)\le\varphi(t^{q-c})=o(1),
\end{eqnarray*}
which proves the localization of $u_{1,1}$ and $u_{1,2}$.

To prove that $u_{1,3}$ is negligible, observe that the contributing
paths do not visit~$Z_t^{({1})}$ and are longer than $|Z_t^{({2})}|(1+t^{-\delta/2})$
(the latter is
true as they pass
through $Z_t^{({2})}$ and leave $\Gamma_t^{({2})}$). Thus, they do not
belong to the set ${A_{1}}$
(defined at the beginning of Section \ref{sec-U12345})
and so, using Lemma \ref{25}, we obtain
\begin{eqnarray*}
\sum_{z\in\mathbb Z^d}u_{1,3}(t,z)\le
U_2(t)+U_3(t)+U_4(t)+U_5(t)=U(t) o(1).
\end{eqnarray*}

(2) Now we consider the {complementary case} $\xi(Z_t^{({2})})>
\xi
(Z_t^{({1})})-t^{q-c}$. Let us
pick $B=\Gamma_t$ and $\Omega=\{Z_t^{({1})},Z_t^{({2})}\}$. We have
\begin{eqnarray*}
\mathfrak{g}_{\Omega,B}
&=&\min\bigl\{\xi\bigl(Z_t^{({1})}\bigr),\xi\bigl(Z_t^{({2})}\bigr)\bigr\}- \max
_{z\in\Gamma
_t\setminus\{Z_t^{({1})},Z_t^{({2})}\}}\xi(z)\\
&>&\xi(Z_t^{({1})})-t^{q-c}- \max_{z\in\Gamma
_t\setminus\{
Z_t^{({1})},Z_t^{({2})}\}}\xi(z)>t^{q-c/2}-t^{q-c},
\end{eqnarray*}
where we used our assumption on the difference between $\xi
(Z_t^{({1})})$ and
$\xi(Z_t^{({2})})$
and Lemma \ref{mb}(i) with the constant $c/2$. By Lemma \ref{spectral}
we now obtain, as $t\to\infty$,
\begin{eqnarray*}
\frac{\sum_{z\in\Gamma_t\setminus\{Z_t^{({1})},Z_t^{({2})}\}
}u_1(t,z)}{\sum
_{z\in\Gamma_t}u_1(t,z)}\le
\varphi(t^{q-c/2}-t^{q-c})=o(1).
\end{eqnarray*}
Again, the denominator will only increase if we replace it by $U(t)$.
For the numerator, we observe that $u_1(t,z)=0$ for all $z\notin\Gamma
_t$ as the paths
corresponding to $u_1$ do not leave $\Gamma_t$. {This completes the proof.}
\end{pf}

\begin{lemma}\label{lem-u1esti2} Almost surely,
\begin{eqnarray*}
\lim_{t\to\infty}\Biggl[U(t)^{-1}\sum_{z\in\mathbb Z
^d}u_2(t,z)\Biggr]
{\mathbh1}\bigl\{Z_t^{({2})}\in\Gamma_t^{({1})}\bigr\}=0.
\end{eqnarray*}
\end{lemma}

\begin{pf} We further split $u_2$ into the three contributions
$u_2=u_{2,1}+u_{2,2}+u_{2,3}$, where
\begin{eqnarray*}
u_{2,j}(t,z)=\mathbb E_0\biggl[\exp\biggl\{\int_0^t\xi(X_s)\,ds\biggr\}
{\mathbh1}
\{X_t=z\}{\mathbh1}_{C_{2j}}\biggr]
\end{eqnarray*}
with
\begin{eqnarray*}
C_{2j}=\bigl\{\tau_{\{Z_t^{({1})},Z_t^{({2})}\}}> t\mbox{ or }\tau
_{\Gamma
_t^{\mathrm{c}}}\le t\bigr\}
\cap\cases{
(A_1\cup A_2\cup A_3)
\cap\{\tau_{\Gamma_t^{\mathrm{c}}}\le t\}, &\quad $j=1$,\cr
(A_1\cup A_2\cup A_3)
\cap\{\tau_{\Gamma_t^{\mathrm{c}}}> t\}, &\quad $j=2$,\cr
(A_4\cup A_5), &\quad $j=3$,}
\end{eqnarray*}
where we recall the events $A_1,\ldots,A_5$ {defined at the beginning
of} Section \ref{sec-U12345}.
Since $A_1,\ldots,A_5$ are pairwise disjoint and $(\bigcup
_{i=1}^5A_i)^{\mathrm{c}}=\varnothing$, the sets $C_{21}$, $C_{22}$
and~$C_{23}$ are pairwise disjoint as well, and their union is equal to
the set
\[
\bigl\{\tau_{\{Z_t^{({1})},Z_t^{({2})}\}}> t\mbox{ or }\tau_{\Gamma
_t^{\mathrm
{c}}}\le t\bigr\}
\]
appearing in the definition of $u_2(t,z)$. Hence, we indeed have
$u_2=u_{2,1}+u_{2,2}+u_{2,3}$.

We argue on the event $\{Z_t^{({2})}\in\Gamma_t^{({1})}\}$, but
only for $u_{2,1}(t,z)$ this condition will be essential.
Each path contributing to $u_{2,1}$ leaves
$\Gamma_t$ and so passes through some point $z\notin\Gamma
^{({1})}_t\cup\Gamma^{({2})}_t$
according to Lemma \ref{mb}(iv). If the path also passes through~$Z_t^{({i})}$ for $i=1$ or $i=2$
then its length must not be less than $|Z_t^{({i})}|(1+t^{-\delta
/2})$. Hence, by Lemma \ref{25},
\begin{eqnarray*}
\sum_{z\in\mathbb Z^d}u_{2,1}(t,z)
&\le U_2(t)+U_3(t)=U(t)o(1).
\end{eqnarray*}
To bound $u_{2,2}$ we observe that as $\tau_{\Gamma_t^{\mathrm{c}}}>
t$, the {alternative}
$\tau_{\{Z_t^{({1})},Z_t^{({2})}\}}> t$ must be satisfied. Hence
we can use Lemma \ref{25} to get
\begin{eqnarray*}
\sum_{z\in\mathbb Z^d} u_{2,2}(t,z)\le U_3(t)=U(t)o(1).
\end{eqnarray*}
Finally, to bound $u_{2,3}$ we simply use Lemma \ref{25} and obtain
\begin{eqnarray*}
\sum_{z\in\mathbb Z^d}u_{2,3}(t,z)\le U_4(t)+U_5(t)=U(t)o(1),
\end{eqnarray*}
which completes the proof.
\end{pf}

\subsection{Second event: $\Phi_t(Z_t^{({1})})$ is close
to $\Phi_t(Z_t^{({2})})$, but $Z_t^{({1})}$ is far from
$Z_t^{({2})}$}\label{sec-Proof2}

In this section we prove Proposition \ref{two}. Again, we decompose
$u=u_1+u_2$ such that~$u_1$ is localized in $Z_t^{({1})}$ and $Z_t^{({2})}$, and that $u_2$
is negligible.
In order to show that
we further decompose $u_1$ and $u_2$ as
\begin{eqnarray*}
u_1(t,z)=\sum_{j=1}^2u_{1,j}(t,z)\quad \mbox{and}\quad
u_2(t,z)=\sum_{j=1}^4 u_{2,j}(t,z),
\end{eqnarray*}
where the functions $u_{i,j}$ are defined by
\begin{eqnarray*}
u_{i,j}(t,z)&=&\mathbb E_0\biggl[\exp\biggl\{\int_0^t\xi(X_s)\,ds\biggr\}
{\mathbh1}
\{X_t=z\}{\mathbh1}_{C_{ij}}\biggr]
\end{eqnarray*}
with
\begin{eqnarray*}
C_{1j}=\cases{
\bigl\{\tau_{Z_t^{({1})}}\le t,\tau_{[\Gamma_t^{({1})}]^{\mathrm
{c}}}>t\bigr\}, &\quad $j=1$,\cr
\bigl\{\tau_{Z_t^{({1})}}> t,\tau_{Z_t^{({2})}}\le t,\tau_{[\Gamma
_t^{({2})}]^{\mathrm{c}}}>t\bigr\}, &\quad $j=2$}
\end{eqnarray*}
and
\begin{eqnarray*}
C_{2j}=\cases{
(A_1\cup A_2\cup A_3)\cap\bigl\{\tau_{Z_t^{({1})}}\le t,\tau
_{[\Gamma
_t^{({1})}]^{\mathrm{c}}}\le t\bigr\}, &\quad $j=1$,\cr
(A_1\cup A_2\cup A_3)\cap\bigl\{\tau_{Z_t^{({1})}}> t,\tau
_{Z_t^{({2})}}> t\bigr\}
, &\quad $j=2$,\cr
(A_1\cup A_2\cup A_3)\cap\bigl\{\tau_{Z_t^{({1})}}> t,\tau
_{Z_t^{({2})}}\le t,\tau
_{[\Gamma_t^{({2})}]^{\mathrm{c}}}\le t\bigr\}, &\quad $j=3$,\cr
(A_4\cup A_5)\cap(C_{11}\cup C_{12})^{\mathrm{c}}, &\quad $j=4$,}
\end{eqnarray*}
where we again recall the definition of the disjoint sets $A_1,\ldots,A_5$
from Section~\ref{sec-U12345}. It is easy to see that the six sets $C_{11}$,
$C_{12}$, $C_{21}$, $C_{22}$, $C_{23}$ and $C_{24}$ are pairwise
disjoint and exhaustive.

\begin{lemma}
Almost surely,
\begin{eqnarray*}
&&\lim_{t\to\infty}
\Biggl[U(t)^{-1}
\sum_{z\in\mathbb Z^d\setminus\{Z_t^{({1})},Z_t^{({2})}\}}
u_1(t,z)\Biggr]\\
&&\qquad{}\times{\mathbh1}\bigl\{\Phi_t\bigl(Z_t^{({1})}\bigr)-\Phi_t\bigl(Z_t^{({2})}\bigr)< a_t\lambda
_t/2,Z_t^{({2})}\notin\Gamma_t^{({1})}\bigr\}=0.
\end{eqnarray*}
\end{lemma}

\begin{pf}
We argue on the event $\{\Phi_t(Z_t^{({1})})-\Phi_t(Z_t^{({2})})<
a_t\lambda_t/2,Z_t^{({2})}\notin\Gamma
_t^{({1})}\}$.
We now fix $t$ large enough and pick $B=\Gamma_t^{({1})}$, $\Omega
=\{Z_t^{({1})}\}$
to study $u_{1,1}$ and $B=\Gamma_t^{({2})}\setminus\{Z_t^{({1})}\}$,
$\Omega=\{Z_t^{({2})}\}$ to study $u_{1,2}$.
Since $Z_t^{({2})}\notin\Gamma_t^{({1})}$ we have {for the first choice}
\begin{eqnarray*}
\mathfrak{g}_{\Omega,B}=\xi\bigl(Z_t^{({1})}\bigr)- \max_{z\in\Gamma
_t^{({1})}\setminus\{Z_t^{({1})}\}}\xi(z)\ge t^{q-c},
\end{eqnarray*}
using parts (i) and (iv) of Lemma \ref{mb}. For the second
choice, we also obtain
\begin{eqnarray*}
\mathfrak{g}_{\Omega,B}=\xi\bigl(Z_t^{({2})}\bigr)-
\max
_{z\in\Gamma_t^{({2})}\setminus\{Z_t^{({1})},Z_t^{({2})}\}}\xi
(z)\ge t^{q-c},
\end{eqnarray*}
by Lemma \ref{mb}(iii) since the condition $\Phi
_t(Z_t^{({1})})-\Phi_t(Z_t^{({2})})<a_t\lambda
_t/2$ is satisfied. By Lemma \ref{spectral}
and using monotonicity of $\varphi$ we now obtain
\begin{eqnarray*}
\frac{\sum_{z\in\Gamma_t^{({1})}\setminus\{Z_t^{({1})}\}
}u_{1,1}(t,z)}{\sum_{z\in\Gamma_t^{({1})}}u_{1,1}(t,z)}
\le\varphi(t^{q-c}) \quad\mbox{and}\quad
\frac{\sum_{z\in\Gamma_t^{(2)}\setminus\{Z_t^{({1})},Z_t^{({2})}\}
}u_{1,2}(t,z)}{\sum_{z\in\Gamma_t^{({2})}\setminus\{Z_t^{({1})}\}}
u_{1,2}(t,z)}\le\varphi(t^{q-c}).
\end{eqnarray*}
Increasing the denominators to $U(t)$ and taking into account the fact
that $u_{1,1}(t,z)=0$
for all $z\notin\Gamma_t^{({1})}$ and $u_{1,2}(t,z)=0$ for all
$z\notin\Gamma_t^{({2})}\setminus\{Z_t^{({1})}\}$ completes the proof.\vadjust{\goodbreak}
\end{pf}

\begin{lemma} Almost surely,
\begin{eqnarray*}
\lim_{t\to\infty}\Biggl[U(t)^{-1}\sum_{z\in\mathbb
Z^d}u_2(t,z)\Biggr]=0.
\end{eqnarray*}
\end{lemma}

\begin{pf} 
Observe that
\begin{eqnarray*}
&& A_1\cap\bigl[\bigl\{\tau_{Z_t^{({1})}}\le t,\tau_{[\Gamma
_t^{({1})}]^{\mathrm{c}}}\le t\bigr\}\\
&&\qquad\hspace*{5.1pt}{}
\cup\bigl\{\tau_{Z_t^{({1})}}> t,\tau_{Z_t^{({2})}}> t\bigr\}
\cup\bigl\{\tau_{Z_t^{({1})}}> t,\tau_{Z_t^{({2})}}\le t,\tau
_{[\Gamma_t^{({2})}]^{\mathrm{c}}}\le t\bigr\}
\bigr]=\varnothing
\end{eqnarray*}
and therefore the union with $A_1$ can be skipped in the definition of
$C_{21}$, $C_{22}$ and~$C_{23}$. By Lemma \ref{25} we obtain, almost surely,
\begin{eqnarray*}
\sum_{z\in\mathbb Z^d}u_{2,j}(t,z)\le U_2(t)+U_3(t)= U(t) o(1)
\qquad\mbox{for }j=1,2,3.
\end{eqnarray*}
Note that, obviously,
$\sum_{z\in\mathbb Z^d}u_{2,4}(t,z)\le U_4(t)+U_5(t)= U(t) o(1)$
almost surely.
\end{pf}

\subsection{Third event: the difference between $\Phi
_t(Z_t^{({1})})$ and
$\Phi_t(Z_t^{({2})})$ is large}\label{sec-Proof3}

In this section we prove Proposition \ref{three}. Here we decompose
$u=u_1+u_2$ and further $u_2=u_{2,1}+u_{2,2}+u_{2,3}$ where
\begin{eqnarray*}
u_1(t,z)&=&\mathbb E_0\biggl[\exp\biggl\{\int_0^t\xi(X_s)\,ds\biggr\}
{\mathbh1}\{
X_t=z\}{\mathbh1}\bigl\{\tau_{Z_t^{({1})}}\le t,\tau_{[\Gamma
_t^{({1})}]^{\mathrm
{c}}}>t\bigr\}\biggr],\\
u_{2,j}(t,z)&=&\mathbb E_0\biggl[\exp\biggl\{\int_0^t\xi(X_s)\,ds\biggr\}
{\mathbh1}
\{X_t=z\}{\mathbh1}_{C_{2j}}\biggr]
\end{eqnarray*}
with
%
\begin{eqnarray*}
C_{2j}=\cases{
(A_1\cup A_2\cup A_3)\cap\bigl\{\tau_{Z_t^{({1})}}> t\bigr\}, &\quad
$j=1$,\cr
(A_1\cup A_2\cup A_3)\cap\bigl\{\tau_{Z_t^{({1})}}\le t,\tau
_{[\Gamma
_t^{({1})}]^{\mathrm{c}}}\le t\bigr\}, &\quad $j=2$,\cr
(A_4\cup A_5)\cap C_1^{\mathrm{c}}, &\quad $j=3$.}
\end{eqnarray*}
Again, it is easy to see that $u$ is equal to the sum of the functions
$u_1$ and $u_{2,1}$, $u_{2,2}$ and $u_{2,3}$.

\begin{lemma}
\label{loc3}
Almost surely,
\begin{eqnarray*}
\lim_{t\to\infty}
\Biggl[U(t)^{-1}
\sum_{z\in\mathbb Z^d\setminus\{Z_t^{({1})}\}}
u_1(t,z)\Biggr]
{\mathbh1}\bigl\{\Phi_t\bigl(Z_t^{({1})}\bigr)-\Phi_t\bigl(Z_t^{({2})}\bigr)\ge a_t\lambda
_t/2\bigr\}=0.
\end{eqnarray*}
\end{lemma}

\begin{pf} We fix $t$ large enough and argue on the event $\{\Phi
_t(Z_t^{({1})})-\Phi_t(Z_t^{({2})})
\ge a_t\lambda_t/2\}$.
Pick $B=\Gamma_t^{({1})}$, $\Omega=\{Z_t^{({1})}\}$. We have
\begin{eqnarray*}
\mathfrak{g}_{\Omega,B}=\xi\bigl(Z_t^{({1})}\bigr)- \max_{\Gamma
_t^{({1})}\setminus\{Z_t^{({1})}\}}\xi(z)
\ge t^{q-c}
\end{eqnarray*}
by Lemma \ref{mb}(ii) since the condition $\Phi
_t(Z_t^{({1})})-\Phi_t(Z_t^{({2})})\ge a_t\lambda_t/2$
is satisfied. Using Lemma \ref{spectral} we obtain
\begin{eqnarray*}
\frac{\sum_{z\in\Gamma_t^{({1})}\setminus\{Z_t^{({1})}\}
}u_1(t,z)}{\sum_{z\in\Gamma_t^{({1})}}u_1(t,z)}
\le\varphi(t^{q-c})=o(1).
\end{eqnarray*}
Increasing the denominators to $U(t)$ and taking into account the fact
that \mbox{$u_1(t,z)=0$}
for all $z\notin\Gamma_t^{({1})}$ completes the proof.
\end{pf}

\begin{lemma}
\label{neg3}
Almost surely,
\begin{eqnarray*}
\lim_{t\to\infty}\Biggl[U(t)^{-1}\sum_{z\in\mathbb
Z^d}u_2(t,z)\Biggr]{\mathbh1}
\bigl\{\Phi_t\bigl(Z_t^{({1})}\bigr)-\Phi_t\bigl(Z_t^{({2})}\bigr)\ge a_t\lambda_t/2\bigr\}=0.
\end{eqnarray*}
\end{lemma}

\begin{pf} We argue on the event $\{\Phi_t(Z_t^{({1})})-\Phi
_t(Z_t^{({2})})\ge a_t\lambda_t/2\}$.
Denote $h_t(z)=0$ and $H_{{t}}=\{Z_t^{({1})}\}$. By Proposition \ref
{iub} and
by Lemma \ref{l10}(b) we have
\begin{eqnarray*}
\frac{1}{t}\log\Biggl[U(t)^{-1}\sum_{z\in\mathbb
Z^d}u_{2,1}(t,z)\Biggr]
&\le&\frac{1}{t}\log\Biggl[U(t)^{-1}\sum_{z\in\mathbb
Z^d}u_{H,h}(t,z)
\Biggr]\\
&\le&\Phi_t\bigl(Z_t^{({2})}\bigr)-\Phi_t\bigl(Z_t^{({1})}\bigr)+O(t^{q-\delta})\\
&\le&-a_t\lambda_t/2+O(t^{q-\delta})\to-\infty.
\end{eqnarray*}
Further, since
$A_1\cap\{\tau_{Z_t^{({1})}}\le t,\tau_{[\Gamma_t^{({1})}]^{\mathrm
{c}}}\le t\}=\varnothing$
the union with $A_1$ can be skipped in the definition of $C_{22}$.
Then by Lemma \ref{25} we obtain, almost surely,
\begin{eqnarray*}
\sum_{z\in\mathbb Z^d}u_{2,2}(t,z)\le U_2(t)+U_3(t)= U(t)o(1).
\end{eqnarray*}
Obviously, we also have
$\sum_{z\in\mathbb Z^d}u_{2,3}(t,z)\le U_4(t)+U_5(t)=U(t)o(1)$
almost surely.
\end{pf}

\section{One point localization in law and concentration sites}\label{s_weak}

In this section we prove Theorems \ref{main_w} and \ref{main_z}, the
convergence assertions for $u(t,Z_t^{({1})})/U(t)$ in probability and for
$(Z_t^{({1})},Z_t^{({2})})/r_t$ in distribution. This easily follows
from our earlier
almost-sure results, using a point process convergence approach.
Background on point processes
and similar arguments can be found in \cite{HMS08}.%

Consider the Radon measure $\mu(dy)=\frac{\alpha \,dy}{y^{\alpha
+1}}$ on $(0,\infty]$ and, for any $r>0$, the point process on
$\mathbb R^d\times(0,\infty]$ given by
%
\begin{equation}
\label{pp}
\zeta_r=\sum_{z\in\mathbb Z^d}\varepsilon_{(z/r,X_{r,z})},\qquad
\mbox
{where }X_{r,z}=\frac{\xi(z)}{r^{d/\alpha}},
\end{equation}
where we write $\varepsilon_x$ for the Dirac measure in $x$.
Furthermore, for any $t$, consider the point process on $\mathbb
R^d\times
(0,\infty]$ given by
\begin{eqnarray*}
\Pi_t=\sum_{z\in\mathbb Z^d\dvtx\Phi_t(z)>0}\varepsilon_{(z/r_t,
\Phi
_t(z)/a_t)}.
\end{eqnarray*}
%
Finally, define a locally compact Borel set
\begin{eqnarray*}
H=\{(x,y)\in\dot{\mathbb R}^d\times(0,\infty] \dvtx y \ge q|x|/2
\},
\end{eqnarray*}
where $\dot{\mathbb R}^d$ is the one point compactification of
$\mathbb R^d$.


\begin{lemma}
\label{pppp}
For each $t$, $\Pi_t$ is a point process on
\[
\widehat H=\dot{\mathbb R}^{d+1}\setminus\bigl(\bigl(\mathbb
R^d\times(-\infty
,0)\bigr) \cup\{(0,0)\} \bigr).
\]
%
As $t\to\infty$,
$\Pi_t$ converges in law to a Poisson process $\Pi$ on $\widehat H$
with intensity measure
\[
\nu(dx,dy)=dx \otimes\frac{\alpha}{(y+q|x|)^{\alpha
+1}}{\mathbh1}_{\{y>0\}}\,dy.
\]
\end{lemma}

\begin{pf} Our first goal is to write $\Pi_t$ as a suitable
transformation of $\zeta_{r_t}$ on $\widehat H$. Introduce $H'=\dot
{\mathbb R}^{d+1}\setminus\{0\}$ and a transformation $T_t\dvtx H \to H'$
given by
\begin{eqnarray*}
T_t(x,y)= \cases{
\bigl(x,y-q|x|-\delta(t,x,y)\bigr), &\quad if $x\neq\infty\mbox{
and }y\neq\infty$,\cr
\infty,&\quad otherwise.}
\end{eqnarray*}
Here $\delta$ is an error function satisfying $\delta(t,x,y)\to0$ as
$t\to\infty$ uniformly in $(x,y)\in K_n^{\mathrm{c}}$, where
\[
K_n=\{(x,y)\in H\dvtx|y|\ge n\}.
\]
Recalling that $\frac{r_t}{ta_t}=\frac{1}{\log t}$, we see that
\begin{eqnarray*}
\frac{\Phi_t(z)}{a_t}
&=& \biggl[\frac{\xi(z)}{a_t}-\frac{|z|}{ta_t}\log a_t-\frac
{|z|}{ta_t} \log\frac{\xi(z)}{a_t}
+\frac{\eta(z)}{ta_t}\biggr]{\mathbh1}\biggl\{\frac{\xi(z)}{a_t}\ge
[\log
t]^{-1}\frac{|z|}{r_t}\biggr\}\\
&=&\biggl[\frac{\xi(z)}{a_t}-\bigl(q+o(1)\bigr) \biggl|\frac{z}{r_t}\biggr|
-\frac{1}{\log t}\biggl|\frac{z}{r_t}\biggr| \log\frac{\xi
(z)}{a_t}+\frac{\eta(z)}{ta_t}\biggr]\\
&&{}\times
{\mathbh1}\biggl\{\frac{\xi(z)}{a_t}\ge{[\log t]^{-1}}\frac
{|z|}{r_t}
\biggr\}.
\end{eqnarray*}
The same fact also implies that
$\frac{\eta(z)}{ta_t}\le|\frac{z}{r_t}|\frac{\log d}{\log t}$
{ for all $z\in\mathbb Z^d$ and $t>0$.}
Hence, we have
%
\begin{eqnarray}
\label{pipi}
\Pi_t = (\zeta_{ r_t}|_H \circ T_t^{-1}) |_{\widehat{H}}
\qquad\mbox{eventually for all $t$.}
\end{eqnarray}
%
{To show the convergence,} we define the transformation $T \dvtx H\to
H'$ by $T(x,y)=(x,y-q|x|)$ if $x\neq\infty$ and
$y\neq\infty$ and $T(x,y)=\infty$ otherwise. By \cite{HMS08},
Lemma~3.7, $\zeta_r|_H$ is a point process in $H$ converging, as
$r\to\infty$, in law to
a Poisson point process $\zeta|_H$ with intensity measure $\mathrm{
Leb}_d\otimes \mu|_H$,
where $\mathrm{Leb}_d$ denotes the Lebesgue measure on $\mathbb R^d$.
Using \textup{(\ref{pipi})},
it now suffices to show that
\[
\zeta_{ r_t}|_H\circ T_t^{-1}\quad\Longrightarrow\quad\zeta|_H\circ T^{-1}
\]
as the Poisson process on the right has the required intensity by a
straightforward
change of coordinates. This convergence follows from \cite{HMS08},
Lemma 2.5,
provided that the conditions (i)--(iii) stated there are
satisfied, which we now check:

\begin{longlist}
\item $T$ is obviously continuous.
\item For each compact set $K'\subset H'$ there is
an open neighborhood $V'$ of zero such that ${K'\subset H'\setminus V'}$.
Since $T(x,y)\to(0,0)$ as $(x,y)\to(0,0)$ and since \mbox{$T_t\to T$}
uniformly on $K_n^{\mathrm{c}}$,
there exists an open neighborhood $V\subset H$ of zero such that
$T(V)\subset V'$ and
$T_t(V)\subset V'$ for all $t$ large enough. Hence,
for $K=H\setminus V$, we obtain $T^{-1}(K')\subset T^{-1}(H'\setminus
V')\subset K$ and
similarly $T^{-1}_t(K')\subset K$ for all $t$.
\item Recall that $\delta(t,x,y)\to0$ uniformly on $K_n^{\mathrm
{c}}$, and observe that
\begin{eqnarray*}
(\mathrm{Leb}_d \otimes \mu)(K_n)&=&\int_{\mathbb R^d}\,dx\int_{n\vee
(q|x|/2)}^{\infty}\frac{\alpha \,dy}{y^{\alpha+1}}\\
&=&(2/q)^{\alpha}\int_{\mathbb R^d}\frac{dx}{((2n/q)\vee|x|)^{\alpha
}}\to0
\end{eqnarray*}
as $n\to\infty$ as $|x|^{-\alpha}$ is integrable away from zero for
$\alpha>d$.\quad\qed
\end{longlist}
\noqed\end{pf}

\begin{lemma}
\label{jointlaw}
We have
\[
\biggl(\frac{Z_t^{({1})}}{r_t},\frac{Z_t^{({2})}}{r_t},\frac{\Phi
_t(Z_t^{({1})})}{a_t},\frac{\Phi_t(Z_t^{({2})})
}{a_t}\biggr)\quad
\Rightarrow\quad\bigl(X^{({1})}, X^{({2})}, Y^{({1})}, Y^{({2})}\bigr),
\]
where the limit random variable has the density
\begin{eqnarray*}
p(x_1,x_2,y_1,y_2)=\frac{\alpha^2\exp\{-\theta y_2^{d-\alpha}\}}
{(y_1+q|x_1|)^{\alpha+1}(y_2+q|x_2|)^{\alpha+1}} {\mathbh1}\{
y_1\geq y_2\}.
\end{eqnarray*}
\end{lemma}

\begin{pf}
It has been computed in the proof of \cite{HMS08}, Proposition 3.8, that
$\nu(\mathbb R^d\times(y,\infty))=\theta y^{d-\alpha}$ for $y>0$.
For any relative compact set $A\subset\widehat H\times\widehat H $
such that
$\mathrm{Leb}_{2d+2}(\partial A)=0$, we obtain by Lemma \ref{pppp},
\begin{eqnarray*}
&&\operatorname{Prob}\biggl(\biggl(\frac{Z_t^{({1})}}{r_t},\frac
{Z_t^{({2})}}{r_t},\frac{\Phi_t(Z_t^{({1})}) }{a_t},\frac{\Phi
_t(Z_t^{({2})})}{a_t}\biggr)\in A\biggr)\\
&&\qquad=\int_A \operatorname{Prob}\bigl(\Pi_t(dx_1\times dy_1)=\Pi_t(d
x_2\times dy_2)=1,\\
&&\hspace*{72.1pt}\Pi_t\bigl(\mathbb R^d\times(y_1,\infty)\bigr)=\Pi
_t\bigl(\mathbb R
^d\times(y_2,y_1)\bigr)=0\bigr)\\
&&\qquad\to\int_A \operatorname{Prob}\bigl(\Pi(dx_1\times dy_1)=1\bigr)
\operatorname{Prob}
\bigl(\Pi(dx_2\times dy_2)=1\bigr)\\
&&\hspace*{48.1pt}{}\times \operatorname{Prob}\bigl(\Pi\bigl(\mathbb R^d\times
(y_1,\infty)\bigr)=0
\bigr) \operatorname{Prob}\bigl(\Pi\bigl(\mathbb R^d\times(y_2,y_1)\bigr)=0\bigr)\\
&&\qquad=\int_A \nu\bigl(\mathbb R^d\times(y_2,\infty)\bigr) \nu(dx_1,dy_1)\nu
(dx_2,dy_2)\\
&&\qquad= \int_A p(x_1,x_2,y_1,y_2)\,dx_1\,dx_2\,dy_1\,dy_2.
\end{eqnarray*}
It remains to notice that
\begin{eqnarray*}
&&\int_{\mathbb R^d\times\mathbb R^d\times\{(y_1> y_2> 0)\}}
p(x_1,x_2,y_1,y_2)\,dx_1\,dx_2\,dy_1\,dy_2\\
&&\qquad=\operatorname{Prob}\bigl(\Pi\bigl(\mathbb R^d\times(0,\infty)\bigr)\ge2\bigr)=1
\end{eqnarray*}
since $\Pi(\mathbb R^d\times(0,\infty))=\infty$ with probability one.
\end{pf}

\begin{pf*}{Proof of Theorem \protect\ref{main_w}}
We use the same
decomposition $u(t,z)=u_1(t, z)+u_2(t,z)$ as we used
to prove Proposition \ref{three}. By Lemmas \ref{loc3} and \ref
{neg3} it suffices to show that
%
\begin{eqnarray}
\label{www}
\lim_{t\to\infty}\operatorname{Prob}\bigl(\Phi
_t\bigl(Z_t^{({1})}\bigr)-\Phi_t\bigl(Z_t^{({2})}\bigr)\ge a_t\lambda_t/2\bigr)=1.
\end{eqnarray}
{Since, by Lemma \ref{jointlaw}, $(\Phi_t(Z_t^{({1})})/a_t,\Phi
_t(Z_t^{({2})})/a_t)$ converges
weakly to a
random variable $(Y^{({1})}, Y^{({2})})$ with density, we
obtain \textup{(\ref{www})} because
$\lambda_t\to0$.}
\end{pf*}

\begin{pf*}{Proof of Theorem \protect\ref{main_z}}
The result follows from
Lemma \ref{jointlaw}
by integrating the density function $p(x_1,x_2,y_1,y_2)$ over all
possible values of $y_1$ and $y_2$. We obtain
\begin{eqnarray*}
p(x_1,x_2)&=& \int_{\{y_1>y_2>0\}}
\frac{\alpha^2\exp\{-\theta y_2^{d-\alpha}\}\,dy_1\,dy_2}
{(y_1+q|x_1|)^{\alpha+1}(y_2+q|x_2|)^{\alpha+1}}\\
&=&\int_0^{\infty} \frac{\alpha\exp\{-\theta y^{d-\alpha}\}\,dy}
{(y+q|x_1|)^{\alpha}(y+q|x_2|)^{\alpha+1}}.
\end{eqnarray*}
This completes the proof of Theorem \ref{main_z}.
\end{pf*}

\section*{Acknowledgments}
We gratefully acknowledge support by the DFG\break \mbox{Forschergruppe 718}
\textit{Analysis and stochastics in complex physical
systems}.

%

%
\printaddresses

\end{document}